\documentclass[a4paper,review]{cas-sc}

\usepackage[section]{placeins} 
\usepackage[capitalise]{cleveref}
\usepackage{booktabs}
\usepackage{xfrac}
\usepackage[sort&compress,square,numbers]{natbib}
\usepackage{graphicx}
\usepackage{graphics}
\usepackage{wrapfig}
\usepackage{float}
\usepackage{subfig}
\usepackage[percent]{overpic}
\usepackage{varwidth}
\usepackage{tikz}
\usepackage{pgfplots}
\usepackage{adjustbox}
\usetikzlibrary{arrows,matrix,positioning,fit}
\usetikzlibrary{shapes,positioning}
\usetikzlibrary{backgrounds}
\usepackage{tikz-layers}
\usepgfplotslibrary{fillbetween}
\usetikzlibrary{intersections}
\pgfplotsset{compat=1.14}
\usepgfplotslibrary{colorbrewer}
\usepgfplotslibrary{patchplots}
\usepgfplotslibrary[colorbrewer]
\usetikzlibrary{pgfplots.colorbrewer}
\usetikzlibrary[pgfplots.colorbrewer]
\usepgfplotslibrary{units}
\usetikzlibrary{spy}
\usepackage{pgfplotstable}
\usepackage{arrayjobx}
\graphicspath{ {./figs/} }
\usetikzlibrary{external}

\newtheorem{theorem}{Theorem}[section]
\newtheorem{corollary}{Corollary}[section]
\newtheorem{lemma}{Lemma}[section]
\newtheorem{definition}{Definition}[section]
\newtheorem{remark}{Remark}[section]

\newcommand\pxvar[2]{\partial_{#2} #1}

\newcommand{\half}{{\frac{1}{2}}}

\begin{document}

\title[mode=title]{Bounds Preserving Temporal Integration Methods for Hyperbolic Conservation Laws}
\shorttitle{Bounds Preserving Temporal Integration Methods}
\shortauthors{T. Dzanic \textit{et al.}}

\author[1]{T. Dzanic}[orcid=0000-0003-3791-1134]
\cormark[1]
\cortext[cor1]{Corresponding author}
\ead{tdzanic@tamu.edu}
\author[2,3]{W. Trojak}[orcid=0000-0002-4407-8956]
\author[1]{F. D. Witherden}[orcid=0000-0003-2343-412X]

\address[1]{Department of Ocean Engineering, Texas A\&M University, College Station, TX 77843}
\address[2]{Department of Aeronautics, Imperial College London, South Kensington, London, SW7 2AZ}
\address[3]{IBM Research, The Hartree Centre, Daresbury, WA4 4AD}

\begin{abstract}
In this work, we present a modification of explicit Runge--Kutta temporal integration schemes that guarantees the preservation of any locally-defined quasiconvex set of bounds for the solution. These schemes operate on the basis of a bijective mapping between an admissible set of solutions and the real domain to strictly enforce bounds. Within this framework, we show that it is possible to recover a wide range of methods independently of the spatial discretization, including positivity preserving, discrete maximum principle satisfying, entropy dissipative, and invariant domain preserving schemes. Furthermore, these schemes are proven to recover the order of accuracy of the underlying Runge--Kutta method upon which they are built. The additional computational cost is the evaluation of two nonlinear mappings which generally have closed-form solutions. We show the utility of this approach in numerical experiments using a pseudospectral spatial discretization without any explicit shock capturing schemes for nonlinear hyperbolic problems with discontinuities. 
\end{abstract}




\maketitle


\section{Introduction}\label{sec:introduction}
This work pertains to the approximation of hyperbolic conservation laws of the form
\begin{equation}\label{eq:gen_hype}
     \begin{cases}
        \pxvar{\mathbf{u}}{t} + \mathbf{\nabla}\cdot\mathbf{F}(\mathbf{u}) = 0, \quad \mathrm{for}\ (\mathbf{x}, t)\in\Omega\times\mathbb{R}_+, \\
        \mathbf{u}(\mathbf{x}, 0)= \mathbf{u}_0(\mathbf{x}), \quad \mathrm{for}\ \mathbf{x}\in \Omega,
    \end{cases}
\end{equation}
where $d$ is some arbitrary space dimension, $\mathbf{u}\in\mathbb{R}^m$ is the solution, $\mathbf{F}(\mathbf{u})\in (\mathbb{R}^m)^d$ is the flux, and $\Omega\subset\mathbb{R}^d$ is the domain. The domain is assumed to be periodic to simplify analysis with respect to the boundary conditions. We assume that there exists an admissible set of solutions to \cref{eq:gen_hype}, and without giving a precise meaning to an admissible solution, which in its own right may be an open problem, we assume that there exists a well-defined set of bounds to the solution which must be satisfied to meet some criteria of admissibility. Furthermore, from \citet{Hoff1985, Frid2001, Lax1954} and related works, we assume that for general nonlinear hyperbolic systems, the notion of admissibility of these sets of solutions implies, to some extent, their convexity. As such, this motivates the development and analysis of numerical schemes in the context of their ability to satisfy convex constraints on the solution. 

The literature on spatial discretizations that enforce some criteria upon the solution is vast, spanning many decades and discretization techniques \citep{Shu2018}. However, a drawback in many of these techniques is their lack of generalizability across the various classes of spatial discretizations. By instead utilizing the method of lines approach \citep{Jameson1981}, modifications to Runge--Kutta (RK) temporal integration schemes have been employed to enforce desirable criteria independently of the spatial discretization \citep{Calvo2006,Iserles2000,Hairer2006}. For more complex criteria, these approaches generally rely on some sort of projection methods (see \citet{Hairer2006} Sec. IV.4) in which the solution is projected onto a desired manifold, with various approaches effectively differing in their choices of the search direction \citep{Hairer2006,Calvo2006}. More recently, relaxation RK methods were introduced in \citet{Ketcheson2019} and shown to preserve any inner product norm by scaling the weights of the underlying RK scheme. This was extended to general convex functionals in \citet{Ranocha2020} and applied to the Euler and Navier--Stokes equations with success.

In contrast to projection and relaxation methods, the objective of this work is to instead introduce a novel approach for explicit RK temporal integration schemes that guarantees the preservation of any locally-defined quasiconvex set of bounds for the solution. The results of the proposed approach differ from projection-type methods in that it does not enforce constraint equalities but instead ensures that the solution only remains bounded (i.e., constraint inequalities). This approach can be considered as a type of nonlinear penalty method for the temporal scheme, which allows for the enforcement of more general constraints and the potential for more efficient numerical implementations in comparison to relaxation methods and incremental direction techniques. Furthermore, the proposed approach has the advantage of being essentially independent of the spatial discretization which can guarantee properties such as strong stability preservation in scenarios where the time step restrictions are not known for the given spatial discretization, albeit with potentially more restrictive conditions on the time step. 

The underlying mechanisms of this approach are conceptually similar to the change-of-variable methods of \citet{Ilinca1998} and \citet{Luo2003} for enforcing positivity of turbulence variables by transforming them to their logarithmic/exponential form. Similarly, the proposed bounds preserving RK (BP-RK) schemes utilize a bijective mapping to transform the solution to an auxiliary space prior to temporal integration, after which the inverse mapping is formed such as to guarantee the resulting solution remains within the bounds. A mass correction step is then performed afterwards to enforce conservation. The resulting temporal schemes are explicit, can be modified to preserve any linear invariant of the system, and recover the order of accuracy of the underlying RK schemes upon which they are built. While the applications of the proposed schemes are shown for hyperbolic conservation laws utilizing RK temporal integration, the general techniques are broadly applicable to a wider range of ordinary and partial differential equations and temporal schemes. 

The remainder of this paper is organized as follows. \Cref{sec:preliminaries} presents the formulations of an abstract spatial discretization and the underlying RK methods. The BP-RK schemes are introduced in \Cref{sec:temporal}, and examples of formulations of bounds are presented in \Cref{sec:bounds}. The proposed schemes are implemented and utilized on a variety of nonlinear hyperbolic systems, with implementation details given in \Cref{sec:implementation} and results shown in \Cref{sec:results}. Conclusions are then drawn in \Cref{sec:conclusion}.
\section{Discretization}\label{sec:preliminaries}

Let $\mathbf{U}_h(t) := \sum_{i \in V} \mathbf{u}_i(t) \phi_i (\mathbf{x})$ be a discrete approximation of the solution $\mathbf{u}$ via some basis $\{\phi\}_{i \in V}$ of a finite-dimensional vector space $X_h$. We consider an explicit semidiscretization of \cref{eq:gen_hype} by an abstract numerical scheme given in the form of 
\begin{equation}\label{eq:gen_disc}
    \pxvar{\mathbf{u}_i}{t} \approx -\sum_{j\in\mathcal{I}(i)} \mathbf{c}_{ij}\cdot\mathbf{F}\big(\mathbf{u}_j \big) = \mathbf{L}\big(\mathbf{u}_i, t\big)
\end{equation}
for $i \in V$, where $\mathcal{I}(i) \subseteq V$ denotes the stencil at $i$ and $\mathbf{c}_{ij}$ is some $\mathbb{R}^d$-valued matrix dependent on the spatial discretization. 

Furthermore, we consider a general explicit Runge--Kutta (RK) method of $s$ stages represented through its Butcher tableau as
\begin{equation}
\renewcommand{\arraystretch}{1.5} 
\begin{array}
{c|c}
c & A\\
\hline
  & b^T 
\end{array},
\end{equation}
where $A \in \mathbb{R}^{s \times s}$ is a strictly lower-triangular matrix and $b,c \in \mathbb{R}^s$. The temporal discretization is given by
\begin{subequations}
    \begin{align}
        \mathbf{u}^{n+1}_i &= \mathbf{u}^n_i + \Delta t \sum_{k=1}^s b_k \mathbf{L}_{ik},\\
        \mathbf{u}_{ij}^* &= \mathbf{u}^n_i + \Delta t \sum_{k=1}^s A_{jk} \mathbf{L}_{ik}, \quad j \in \{1, \hdots, s\},
    \end{align}\label{eq:rk}
\end{subequations}
where $\mathbf{L}_{ik} = \mathbf{L}(\mathbf{u}_{ik}^*, t^n + c_k\Delta t)$, $\mathbf{u}_i^n \approx \mathbf{u}_i(t^n)$, $\mathbf{u}_i^{n+1} \approx \mathbf{u}_i(t^{n+1})$, and $t^{n+1} = t^n + \Delta t$ for some time step $\Delta t > 0$ and $n \in \mathbb{N}$.

\section{Bounds Preserving Temporal Integration}\label{sec:temporal}

Let $\mathcal{B} \subseteq \mathbb{R}^m$ be some open, quasiconvex, non-empty set of admissible solutions to \cref{eq:gen_hype}. More generally, let there exist a unique admissible set $\mathcal{B}_i$ for each $i \in V$. We state that the temporal integration scheme is bounds preserving if for any solution $\mathbf{u}_i^n \in \mathcal{B}_i$, there exists a sufficiently small time step $\Delta t > 0$ such that $\mathbf{u}_i^{n+1} \in \mathcal{B}_i\ \forall\ i \in V$ irrespective of the spatial discretization $\mathbf{L} (\mathbf{u}_i, t^n)$. It is clear that in their general form, the explicit RK schemes given by \cref{eq:rk} are not guaranteed to be bounds preserving. 

To address this, consider the mapping $\mathbf{G}_i : \mathcal{B}_i \mapsto \mathbb{R}^m$ for some $i \in V$. We define an auxiliary variable $\mathbf{w} \in \mathbb R^m$ such that $\mathbf{G}_i$ is bijective with respect to $\mathbf{w}$, which yields the relations
\begin{equation}\label{eq:mapping}
    \mathbf{w} := \mathbf{G}_i(\mathbf{u}), \quad \quad \mathbf{u} = \mathbf{G}_i^{-1}(\mathbf{w}),
\end{equation}
for some $\mathbf{u} \in \mathcal B_i$. If we further assume that $\mathbf{G}_i\in C^1(\mathcal{B}_i)$, then an auxiliary semidiscrete equation can be given as
\begin{equation}\label{eq:w_system}
    \pxvar{\mathbf{w}_i}{t} = \mathbf{G}_i'({\mathbf{u}_i}) \mathbf{L}\big(\mathbf{u}_i, t\big),
\end{equation}
where $\mathbf{G}_i'({\mathbf{u}_i})$ denotes the Jacobian of the mapping with respect to $\mathbf{u}$. If the Jacobian is bounded, the auxiliary system is exactly the image of 
\cref{eq:gen_disc}. Since the set $\mathcal{B}_i$ is open, given a solution $\mathbf{u}_i$ at some time $t^*$, the Jacobian is guaranteed to be bounded over the interval $[t^*, t^* + \Delta t]$ in the limit as $\Delta t \to 0$. Utilizing this auxiliary form, we introduce an intermediate temporal update as
\begin{equation}\label{eq:rk_w}
    \overline{\mathbf{u}}^{n+1}_i = \mathbf{G}_i^{-1}(\mathbf{w}^{n+1}_i) = \mathbf{G}_i^{-1}\bigg [ \mathbf{w}^n_i + \Delta t \sum_{k=1}^s b_k \mathbf{G}_i'({\overline{\mathbf{u}}_{ik}^*}) \overline{\mathbf{L}}_{ik}\bigg ],
\end{equation}
where 
\begin{equation}\label{eq:rk_w2}
    \overline{\mathbf{u}}_{ij}^* =  \mathbf{G}_i^{-1}\bigg [ \mathbf{w}^n_i + \Delta t \sum_{k=1}^s A_{jk} \mathbf{G}_i'({\overline{\mathbf{u}}_{ik}^*}) \overline{\mathbf{L}}_{ik}\bigg ],
\end{equation}
and
\begin{equation}
    \overline{\mathbf{L}}_{ik} = \mathbf{L}(\overline{\mathbf{u}}_{ik}^*, t^n + c_k\Delta t).
\end{equation}
These intermediate states utilize the property that the range of $\mathbf{G}_i^{-1}(\mathbf{w})$ is $\mathcal B_i$.
\begin{lemma}[Bounds Preservation of the Intermediate States]\label{lem:bounds}
    For any $i \in V$, let $\mathcal{B}_i \subseteq \mathbb{R}^m$ be some open, quasiconvex, non-empty set and let $\mathbf{G}_i : \mathcal{B}_i \mapsto \mathbb{R}^m$ be a bijective $C^1(\mathcal{B}_i)$ mapping. Given a solution $\mathbf{u}^{n}_i \in \mathcal B_i$, there exists a finite time step $\Delta t > 0$ such that $\overline{\mathbf{u}}^{n+1}_i \in \mathcal B_i$ and $\overline{\mathbf{u}}_{ij}^* \in \mathcal B_i$ for all $j \in \{1, \hdots, s\}$ given the temporal update in \cref{eq:rk_w,eq:rk_w2}.
\end{lemma}
It can be seen that there exists a bijective $C^1(\mathcal{B})$ mapping $\mathbf{G}$ for any open, quasiconvex, non-empty set $\mathcal{B} \subseteq \mathbb{R}^m$. As such, it is possible to construct a mapping that preserves any locally-defined quasiconvex set of bounds for the solution. The construction of these bounds and mappings is further explored in \Cref{sec:bounds}.

For any temporal integration scheme, it is essential for the scheme to at least preserve linear invariants (e.g., total mass in a periodic domain). We present the definition of linear invariant preservation through the notion of an arbitrary linear invariant residual.  

\begin{definition}[Linear Invariant Residual]\label{def:lin_inv_res}
    Let $Q(\mathbf{u})$ be some linear invariant of $\mathbf{u}$ (i.e., $Q(\mathbf{u}) = \mathbf{c}\cdot\mathbf{u}$, $\mathbf{c} \in \mathbb{R}^m \setminus \{ \mathbf{0} \}$). The linear invariant residual is defined as
    \begin{equation}
        R(\mathbf{u}, \mathbf{u}') = Q(\mathbf{u}) - Q(\mathbf{u}') = Q(\mathbf{u} - \mathbf{u}').
    \end{equation}
\end{definition}

\begin{definition}[Linear Invariant Preservation]\label{def:lin_inv_conv}
    Let $m_i$ be some positive quantity dependent on the spatial discretization such that 
    \begin{equation*}
        \sum_{i \in V} m_i \mathbf{u}_i^n \approx \int_\Omega \mathbf{U}_h (t^n).
    \end{equation*}
    The scheme is said to preserve linear invariants for a periodic domain if
    \begin{equation}
        \sum_{i \in V} m_i Q(\mathbf{u}_i^n) = \sum_{i \in V} m_i Q(\mathbf{u}_i^{n'}),
    \end{equation}
    for any linear invariant $Q(\mathbf{u})$ and $n, n' \in \mathbb{N}$, which may be identically expressed as 
    \begin{equation}
        \sum_{i \in V} m_i Q(\mathbf{u}_i^n - \mathbf{u}_i^{n'}) =\sum_{i \in V} m_i R(\mathbf{u}_i^n, \mathbf{u}_i^{n'}) = 0.
    \end{equation}
\end{definition}

Although the proposed intermediate temporal integration scheme is bounds preserving via the intermediate state $\overline{\mathbf{u}}^{n+1}_i$, it can be seen that it does not necessarily preserve any linear invariant of the system due to the nonlinearity of the mappings. Therefore, we define a linear invariant preserving temporal update as
\begin{equation}\label{eq:time_scheme}
    \mathbf{u}^{n+1}_i = \overline{\mathbf{u}}^{n+1}_i + \mathbf{S}_i,
\end{equation}
where $\mathbf{S}_i \in \mathbb R^m$ is an additional term to account for the mass defect in the intermediate states that is of similar form to the corrections performed in \citet{Kuzmin2000}.
The global mass defect is defined as 
\begin{equation}\label{eq:sbar}
    \overline{\mathbf{S}} := \sum_{i \in V} m_i (\mathbf{u}_i^n - \overline{\mathbf{u}}_i^{n+1}),
\end{equation}
from which a component-wise unit vector can be given as 
\begin{equation}\label{eq:n}
    \mathbf{n} := \overline{\mathbf{S}}/\| \overline{\mathbf{S}} \|_2.
\end{equation}
If $\mathbf{S}_i$ is along any arbitrary unit vector, there exists a maximum vector length which can be supported such that $\overline{\mathbf{u}}^{n+1}_i + \mathbf{S}_i \in \mathcal B_i$. This length is formally defined by the following.
\begin{definition}[Set Distance]\label{def:dist}
    For a set $X\subset \mathbb R^m$, the Euclidean distance from a state $\mathbf{u}\in X$ to the boundary of the set, $\partial X$, along some unit vector $\mathbf{n}\in\mathbb{B}^{m-1}(\mathbf{0},1)$ is defined as 
    \begin{equation}
        D_X(\mathbf{u}, \mathbf{n}) = \underset{\gamma \geq 0}{\mathrm{arg\ min}} \Big (\inf_{\mathbf{x}\in\partial X}\|\mathbf{u} + \gamma \mathbf{n} - \mathbf{x}\|_2 \Big).
    \end{equation}
\end{definition}
Along the $\mathbf{n}$ direction, this maximum length is defined as
\begin{equation}\label{eq:abij}
    \gamma_{i}^* :=  D_{\mathcal B_i} (\mathbf{\overline{u}}^{n+1}_i, \mathbf{n}).
\end{equation}
By setting $\mathbf{S}_i$ as 
\begin{equation}\label{eq:si}
    \mathbf{S}_i = \gamma_i\mathbf{n},
\end{equation}
it can be seen from \cref{def:dist} that for any $\gamma_i \in [0,\gamma_{i}^*]$, $\mathbf{u}^{n+1}_i =\overline{\mathbf{u}}^{n+1}_i + \mathbf{S}_i\in \mathcal B_i$. 
We therefore set $\gamma_i$ as 
\begin{equation}\label{eq:gammai}
    \gamma_i = 
      \gamma_{i}^* \| \overline{\mathbf{S}} \|_2/ \sum_{j \in V} m_j \gamma_{j}^*,
\end{equation}
as this will be shown in \cref{thm:bounds,thm:lin_inv} to guarantee a bounds preserving (i.e., $\gamma_i \leq \gamma_i^*$) mass correction step. 

The overall temporal update utilizing the proposed schemes can then be summarized by the following steps:
\begin{enumerate}
    \item For some solution $\{\mathbf{u}_i^{n}\}_{i \in V}$, form a set of bounds $\mathcal B_i$ (which may be unique for each $i \in V$) such that $\mathbf{u}_i^n \in \mathcal B_i \ \forall \ i \in V$.
    \item Create a set of mappings $\mathbf{G}_i : \mathcal{B}_i \mapsto \mathbb{R}^m$.
    \item Transform the solution to the auxiliary space: $\mathbf{w}_i^{n} = \mathbf{G}_i \left (\mathbf{u}_i^{n} \right )$.
    \item Perform the temporal update in auxiliary space and revert the transformation as per \cref{eq:rk_w} to recover $\overline{\mathbf{u}}_i^{n+1}$.
    \item Evaluate the mass defect as per \cref{eq:sbar}.
    \item Apply corrections to preserve linear invariants as per \cref{eq:time_scheme,eq:si,eq:gammai} to recover $\mathbf{u}_i^{n+1}$.
\end{enumerate}

With this formulation, we now move on to state and prove the properties of the proposed schemes. The subsequent theorems utilize the following assumptions:

\begin{enumerate}\label{item:assumptions}
    \item There exists an open, quasiconvex, non-empty set $\mathcal{B}_i \subseteq \mathbb{R}^m  \ \forall \     i \in V$ such that $\mathbf{u}_i^n \in \mathcal B_i \ \forall \ i \in V$.
    \item There exists a bijective $C^\infty(\mathcal{B}_i)$ mapping,
        $\mathbf{G}_i: \mathcal{B}_i \mapsto \mathbb{R}^m$, such that for all $i \in V$, the intermediate temporal update defined by \cref{eq:rk_w} results in $\overline{\mathbf{u}}_i^{n+1} \in \mathcal B_i$ if ${\mathbf{u}}_i^{n} \in \mathcal B_i$.
    \item Periodic boundary conditions are enforced (i.e., $\Omega$ is a $d$-torus).
    \item The underlying spatial discretization $\mathbf{L}(\mathbf{u}, t)$ preserves linear invariants (i.e., $\sum_{i \in V} m_i \mathbf{R}(\mathbf{u}_i^n, \mathbf{u}_i^{n'}) = \mathbf{0}$ for any $n, n'$).
\end{enumerate}

\begin{theorem} [Convergence] \label{thm:convergence}
The temporal integration scheme defined by \cref{eq:rk_w,eq:rk_w2,eq:time_scheme} converges at the rate of the base RK scheme defined in \cref{eq:rk}.
\end{theorem}
\newproof{pot_conv}{Proof of \cref{thm:convergence}}
\begin{pot_conv}
Let $p \geq 1$ be the order of convergence of the base RK scheme defined in \cref{eq:rk} such that for the auxiliary system defined by \cref{eq:w_system}, the relation
\begin{equation}
    \mathbf{w}^{n+1}_i - \mathbf{w}^{n}_i = \Delta t \sum_{k=1}^s b_k \mathbf{T}_k =  \int_{t^n}^{t^{n+1}}\partial_t \mathbf{w}_i(\tau)\ \mathrm{d}\tau + \mathcal O (\Delta t^{p+1}),
\end{equation}
holds for any $i \in V$, where 
\begin{equation*}
    \mathbf{T}_k = \mathbf{G}_i'({\overline{\mathbf{u}}_{ik}^*}) \overline{\mathbf{L}}_{ik}.
\end{equation*} 
The error estimate for the temporal integration scheme defined by \cref{eq:rk_w,eq:rk_w2,eq:time_scheme} can be given as 
\begin{equation}\label{eq:conv_error}
\mathbf{u}^{n+1} - \mathbf{u}^{n} =  \mathbf{G}^{-1} \big (\mathbf{w}^{n} +  \Delta t \sum_{k=1}^s b_k \mathbf{T}_k \big)  - \mathbf{G}^{-1} \big (\mathbf{w}^{n}\big )  + \mathbf{S}.
\end{equation}
Note that the subscript is dropped for brevity. The Taylor series of $\mathbf{H} (\mathbf{w}) = \mathbf{G}^{-1}(\mathbf{w})$ can be expanded around $\mathbf{w} = \mathbf{w}^n$ and evaluated at $\mathbf{w}^{n+1} = \mathbf{w}^n +  \Delta t \sum_{k=1}^s b_k \mathbf{T}_k$ to yield
\begin{equation}
    \mathbf{G}^{-1} \big (\mathbf{w}^{n} + \Delta t \sum_{k=1}^s b_k \mathbf{T}_k\big) = \mathbf{H} \big (\mathbf{w}^{n} \big ) + \sum_{k=1}^{\infty} \frac{1}{k!} \mathbf{H}^{(k)}\big (\mathbf{w}^{n} \big ) \bigg( \Delta t \sum_{j=1}^s b_j \mathbf{T}_j \bigg)^k.
\end{equation}
From assumption 2, the higher-order terms are well-defined as $\mathbf{G}(\mathbf{u}) \in C^\infty(\mathcal{B})$. This may be substituted into \cref{eq:conv_error} to give
\begin{equation}\label{eq:error}
    \mathbf{u}^{n+1} - \mathbf{u}^{n} = \sum_{k=1}^{\infty} \frac{1}{k!} \mathbf{H}^{(k)}\big (\mathbf{w}^{n} \big ) \bigg( \Delta t \sum_{j=1}^s b_j \mathbf{T}_j  \bigg)^k  + \mathbf{S}.
\end{equation}
From \citet{Ranocha2020} (Theorem 2.12 with $\psi = \mathbf{H}(\mathbf{w}, t)$), this expansion is accurate to the equivalent order due to the required accuracy of the underlying RK method as a quadrature rule.
\begin{equation}
    \sum_{k=1}^{\infty} \frac{1}{k!} \mathbf{H}^{(k)}\big (\mathbf{w}^{n} \big ) \bigg( \Delta t \sum_{j=1}^s b_j \mathbf{T}_j \bigg)^k = \int_{t^n}^{t^{n+1}}\partial_t \mathbf{u}_i(\tau)\ \mathrm{d}\tau + \mathcal O (\Delta t^{p+1}).
\end{equation}

It now remains to be shown that $\mathbf{S} \approx \mathcal O (\Delta t^{p+1})$. For some arbitrary $i \in V$, this term may be expressed as
\begin{equation}
    \mathbf{S} = \mathbf{S}_i = \zeta_i \sum_{j \in V} m_j (\mathbf{u}_i^n - \overline{\mathbf{u}}_i^{n+1}) = \zeta_i \sum_{j \in V} m_j \big (\mathbf{u}_i^n - \mathbf{G}^{-1}(\mathbf{w}^{n+1}_i ) \big), 
\end{equation}
where $\zeta_i = \gamma_{i}^* / \sum_{j \in V} m_j \gamma_{j}^*$. Utilizing the fact that $m_j$ and $\zeta_i$ are independent of $\Delta t$ and that, from \cref{eq:mapping}, $-(\mathbf{u}_i^n - \mathbf{G}^{-1}(\mathbf{w}^{n+1}_i ))$ is identical to \cref{eq:error} with $\mathbf{S}_i = 0$, it follows that 
\begin{multline}
    \mathbf{S}_i =  -\zeta_i \sum_{j \in V} m_j \bigg [ \int_{t^n}^{t^{n+1}}\partial_t \mathbf{u}_j(\tau)\ \mathrm{d}\tau   + \mathcal O(\Delta t^{p+1})\bigg ] \\= -\zeta_i  \int_{t^n}^{t^{n+1}} \sum_{j \in V} m_j \partial_t \mathbf{u}_j(\tau)\ \mathrm{d}\tau  + \mathcal O(\Delta t^{p+1}) = \mathcal O(\Delta t^{p+1})
\end{multline}
under the assumption that $\mathbf{L}(\mathbf{u}_j, t)$ is a consistent approximation of $\partial_t \mathbf{u}_j(t)$. As a result, we obtain

\begin{equation}
    \mathbf{u}^{n+1} - \mathbf{u}^{n} = \int_{t^n}^{t^{n+1}}\partial_t \mathbf{u}(\tau)\ \mathrm{d}\tau + \mathcal O (\Delta t^{p+1}).
\end{equation}
\end{pot_conv}

\begin{corollary}[Consistency]\label{thm:consistency}
From \cref{thm:convergence}, the temporal integration scheme defined by \cref{eq:rk_w,eq:rk_w2,eq:time_scheme} is consistent in the sense that
\begin{equation*}
    \underset{\Delta t \to 0}{\mathrm{lim}} \frac{\mathbf{u}_i^{n+1} - \mathbf{u}_i^{n}}{\Delta t} = \frac{\partial \mathbf{u}_i^n}{\partial t}
\end{equation*}
for any $i \in V$.
\end{corollary}

\begin{theorem}[Linear Invariant Preservation]\label{thm:lin_inv}
 The temporal scheme defined by \cref{eq:rk_w,eq:rk_w2,eq:time_scheme} preserves any linear invariant $Q(\mathbf{u})$.
\end{theorem}

\newproof{pot_lip}{Proof of \cref{thm:lin_inv}}
\begin{pot_lip}
For global linear invariant preservation, it is necessary that 
  \begin{equation*}
      \sum_{i \in V} m_i R(\mathbf{u}_i^{n+1}, \mathbf{u}_i^n) = \sum_{i \in V} m_i Q(\mathbf{u}_i^{n+1}- \mathbf{u}_i^n) =
      0,
  \end{equation*}
  which can be equivalently expressed as 
  \begin{equation*}
        Q(\overline{\mathbf{S}}) - \sum_{i \in V} m_i Q(\mathbf{S}_i) = 0.
  \end{equation*}
  The conclusion follows readily from the definition of $\mathbf{S}_i$.
  \begin{equation}
      \sum_{i \in V} m_i Q(\mathbf{S}_i) =  Q \bigg (\sum_{i \in V} m_i \gamma_{i}^* \overline{\mathbf{S}}/ \sum_{j \in V} m_j \gamma_{j}^* \bigg) = Q(\overline{\mathbf{S}}).
  \end{equation}
\end{pot_lip}

\begin{remark} [Boundary Contributions]\label{rmk:bc}
The proof of \cref{thm:lin_inv} is contingent on the assumption that periodic boundary conditions are enforced. It is trivial to extend the definition of $\mathbf{S}_i$ and the subsequent proof of linear invariant preservation to non-periodic domains with boundary contributions. However, this requires some dependency on the spatial discretization and does place some restriction on the formulation of the bounds as they must be able to support the boundary contributions. 
\end{remark}

\begin{remark} [Local Mass Conservation]
The proposed approach for the linear invariant preserving temporal update is designed to preserve \textit{global} mass on a periodic domain. It is possible to extend the definition of $\mathbf{S}_i$ to enforce a local mass balance on a subdomain through some formulation of the ingoing/outgoing fluxes. In similar vein to \cref{rmk:bc}, this requires some dependency on the spatial discretization and does place some restriction on the formulation of the bounds as they must be able to support the contribution of the flux balance. As a result, the general proofs of bounds preservation and linear invariant preservation would not hold without stricter conditions on the bounds. 
\end{remark}

\begin{theorem}[Bounds Preservation]\label{thm:bounds}
  The temporal scheme defined by \cref{eq:rk_w,eq:rk_w2,eq:time_scheme} is bounds preserving for all $\{ \mathcal{B}_i\}_{i \in V}$.
\end{theorem}

\newproof{pot_bounds}{Proof of \cref{thm:bounds}}
\begin{pot_bounds}
    The temporal scheme defined by \cref{eq:rk_w,eq:rk_w2,eq:time_scheme} is bounds preserving if 
    \begin{equation}
    \mathbf{u}^{n+1}_i = \overline{\mathbf{u}}^{n+1}_i + \mathbf{S}_i \in \mathcal B_i \quad \forall \ i \in V.
    \end{equation}
    From \cref{eq:si,eq:gammai}, this condition is satisfied if
    \begin{equation}
        \mathbf{S}_i = \gamma_{i}^*  \overline{\mathbf{S}}/ \sum_{j \in V} m_j \gamma_{j}^* \leq \gamma_{i}^* \mathbf{n} = \gamma_{i}^* \overline{\mathbf{S}} / \| \overline{\mathbf{S}} \|_2,
    \end{equation}
    which may be expressed as 
    \begin{equation}
        \sum_{j \in V} m_j \gamma_{j}^* \geq \| \overline{\mathbf{S}} \|_2.
    \end{equation}
    Here, we switch the index from $j$ to $i$ for consistency with the general notation. Using $\mathbf{n} = \overline{\mathbf{S}}/ \|\overline{\mathbf{S}}\|_2$ and 
    \begin{equation}
        \sum_{i \in V} m_i \gamma_{i}^* = \Big \| \sum_{i \in V} m_i \gamma_{i}^* \mathbf{n} \Big \|_2,
    \end{equation}
    it is sufficient to show that 
    
    \begin{equation}
        \bigg( \overline{\mathbf{S}} - \sum_{i \in V} m_i \gamma_{i}^* \mathbf{n} \bigg)\cdot \mathbf{n} \leq 0. 
    \end{equation}
    From \cref{eq:sbar}, we get
    \begin{equation}
         \bigg(\overline{\mathbf{S}} - \sum_{i \in V} m_i \gamma_{i}^* \mathbf{n} \bigg)\cdot \mathbf{n} = \bigg(\sum_{i \in V} m_i (\mathbf{u}_i^n - \overline{\mathbf{u}}_i^{n+1} - \gamma_{i}^* \mathbf{n})\bigg)\cdot \mathbf{n},
    \end{equation}
    from which we obtain a sufficient condition
    \begin{equation}
         \mathbf{n}\cdot (\mathbf{u}_i^n - \overline{\mathbf{u}}_i^{n+1}) \leq \gamma_{i}^* \quad \forall \ i \in V.
    \end{equation}
    The left-hand side attains its maximum value when
    \begin{equation}
        \mathbf{u}_i^n - \overline{\mathbf{u}}_i^{n+1} = \mathbf{n} \|\mathbf{u}_i^n - \overline{\mathbf{u}}_i^{n+1}\|_2,
    \end{equation}
    which is equivalent to
    \begin{equation}
        \mathbf{u}_i^n  = \overline{\mathbf{u}}_i^{n+1} + \mathbf{n} \|\mathbf{u}_i^n - \overline{\mathbf{u}}_i^{n+1}\|_2.
    \end{equation}
    Given that $\mathbf{u}_i^n \in \mathcal B_i$, from \cref{def:dist}, it can be seen that
    \begin{equation}
        \|\mathbf{u}_i^n - \overline{\mathbf{u}}_i^{n+1}\|_2 \leq \gamma_i^*.
    \end{equation}
    Therefore, we obtain
    \begin{equation}
        \mathbf{n}\cdot (\mathbf{u}_i^n - \overline{\mathbf{u}}_i^{n+1}) \leq   \|\mathbf{u}_i^n - \overline{\mathbf{u}}_i^{n+1}\|_2 \leq \gamma_{i}^*.
    \end{equation}
\end{pot_bounds}
\section{Formulations for Mappings and Bounds}\label{sec:bounds}
The choice of the mapping function $\mathbf{G} (\mathbf{u})$ dictates the properties of the proposed temporal integration schemes. In this section, we present some general mapping functions for various constraints that are of interest in numerical schemes as well as some examples of bounds for hyperbolic conservation laws. 

\subsection{General Mappings}\label{ssec:mappings}
For an arbitrary one-sided constraint on a scalar variable $u$ such that $u > a$ for any $a \in \mathbb R$, the admissible set is defined as 
\begin{equation}
    \mathcal B(u) := (a, \infty).
\end{equation}
A similar bound for $u < a$ may be formed by considering a sign change in $u$. One such example of a bijective mapping function that satisfies the condition $G^{-1}(w) \in \mathcal B$ is 

\begin{equation}
    G(u) = \mathrm{log} ( u - a), \quad \quad
    G^{-1}(w) = \mathrm{exp} (w) + a.
\end{equation}

For a two-sided constraint (i.e., $a < u < b$, $a,b \in \mathbb R$, $b > a$), the admissible set is defined as
\begin{equation}
    \mathcal B(u) := (a, b),
\end{equation}
for which a mapping function can be given by
\begin{equation}
    G(u) = \mathrm{tanh}^{-1} \bigg ( 2\frac{u-a}{b-a} - 1 \bigg ), \quad \quad
    G^{-1}(w) = \frac{a+b}{2} + \frac{b-a}{2} \mathrm{tanh} (w).
\end{equation}

For a vector-valued solution $\mathbf{u} \in \mathbb R^m$, a constraint that is of particular use for hyperbolic conservation laws is given by
\begin{equation}
    \mathcal B(\mathbf{u}) := \{ \mathbf{u} \ \mid \ \| \mathbf{u} \|_2 < r_0\}
\end{equation}
for some $r_0 > 0$, which enforces the condition that the solution (or some subset thereof) exists in an open $m$-ball of radius $r_0$. This mapping is realized through an intermediary mapping, $\mathbf{F} : \mathcal B \to \mathcal Z$, which maps the $m$-ball to an $m$-cube, after which the mapping $\mathcal Z \to \mathbb R^m$ can be formulated as $m$ independent two-sided scalar constraints. Various approaches exist that satisfy the necessary conditions of the intermediary mapping, from elliptic mapping to Schwarz-Christoffel conformal mapping. An example of a mapping for $m=2$ using an elliptic approach is given as 
\begin{align}
    \mathbf{F} \begin{pmatrix}u_1 \\ u_2 \end{pmatrix} &= 
    \begin{pmatrix} 
    \frac{1}{2} \sqrt{2 + \delta \mathbf{u} + 2 \sqrt{2} \big(\frac{u_1}{r_0}\big )}  -
    \frac{1}{2} \sqrt{2 + \delta \mathbf{u} - 2 \sqrt{2} \big(\frac{u_1}{r_0}\big )} \\
    \frac{1}{2} \sqrt{2 - \delta \mathbf{u} + 2 \sqrt{2} \big(\frac{u_2}{r_0}\big )}  -
    \frac{1}{2} \sqrt{2 - \delta \mathbf{u} - 2 \sqrt{2} \big(\frac{u_2}{r_0}\big )}
    \end{pmatrix}\\
    \mathbf{F}^{-1} \begin{pmatrix}z_1 \\ z_2 \end{pmatrix} &= 
    \begin{pmatrix} 
    r_0 z_1 \sqrt{1 - \frac{1}{2} z_2^2} 
    \\ r_0 z_2 \sqrt{1 - \frac{1}{2} z_1 ^2} 
    \end{pmatrix} 
\end{align}
where $\delta \mathbf{u} = \Big(\frac{u_1}{r_0}\Big )^2 - \Big(\frac{u_2}{r_0}\Big)^2$.

\subsection{Bounds}
Although the proposed schemes allow for arbitrary bounds to be placed on any system in question, the focus of this work is on hyperbolic systems of conservation laws. Judicious choices of the bounds allow for these schemes to enforce certain key properties of hyperbolic conservation laws (essentially) independently of the spatial discretization. We briefly present some examples of these bounds.  


\subsubsection{Positivity Preserving}
A common constraint in hyperbolic systems is that $\mathbf{u}$ (or some subset thereof) remains strictly positive across the entire domain (e.g., density and total energy in the Euler equations, water height in the shallow water equations, etc.). Schemes that strictly enforce this condition, usually via the spatial discretization, are generally referred to as positivity preserving. This condition can be enforced through the bounds
\begin{equation}
    \mathcal B_i(\mathbf{u}) := (0, \infty)^m \quad \forall \ i \in V.
\end{equation}

\subsubsection{Discrete Maximum Principle Satisfying}\label{ssec:dmp}
Enforcing the maximum principle property \citep{Lax1954} at the discrete level is a feature that is beneficial to the application and analysis of numerical schemes as it is necessary for the unique entropy solution of scalar hyperbolic conservation laws to obey this property \citep{Dafermos2010}. The formulation of these bounds is not necessarily unique and does have some dependency on the spatial discretization. For instance, one such definition can be 
\begin{equation}
    \mathbf{u}_i^{\mathrm{min}} = \underset{j \in \mathcal A (i)}{\mathrm{min}} \mathbf{u}_j, \quad \quad \mathbf{u}_i^{\mathrm{max}} = \underset{j \in \mathcal A (i)}{\mathrm{max}} \mathbf{u}_j,
\end{equation}
where $\mathcal A (i)$ is some domain of influence of point $\mathbf{x}_i$. A possible formulation of $\mathcal A(i)$ is the numerical domain of influence $\mathcal I(i)$ (i.e., the set of indices for which their associated shape function has support on $\mathbf{x}_i$) or some physical domain of influence related to the relative nodal spacing and propagation speed of the system (i.e., the direct Voronoi neighbors of $\mathbf{x}_i$). Stricter definitions for the bounds can draw on the works of \citet{Lax1954,Nessyahu1990,Guermond2019} and utilize the direction of information propagation within the system (see \Cref{ssec:idp}). Regardless of the choice of the formulation, the bounds are then taken as such to enforce the discrete maximum principle.
\begin{equation}
    \mathcal B_i(\mathbf{u}) := (\mathbf{u}_i^{\mathrm{min}}, \mathbf{u}_i^{\mathrm{max}}).
\end{equation}

\subsubsection{Entropy Dissipative}
Let $\eta(\mathbf{u})$ be some convex entropy functional of \cref{eq:gen_hype}. For entropy dissipative systems, if the condition  
\begin{equation*}
    \frac{\partial}{\partial t} \int_{\mathcal D}\eta(\mathbf{u}) \leq 0
\end{equation*}
is not satisfied at least in a discrete sense, it may lead to solutions that are unphysical and qualitatively incorrect. To ensure that this property is not violated, the bounds can be formulated as 
\begin{equation}
    \mathcal B_i(\mathbf{u}) := \{ \mathbf{u} \ \mid \ \eta (\mathbf{u}) < \eta_i^\mathrm{max} \},
\end{equation}
where $\eta_i^\mathrm{max} \in \mathbb R$ is some maximum entropy. A discrete local entropy inequality can by enforced through a discrete maximum principle bound on the entropy functional.
To make the computation of the mapping function and its associated metrics tractable, it is beneficial for these terms to be algebraically defined to circumvent the need for iterative methods. For many convex functionals, identical constraints can be enforced through alternate formulations of the functional that are more amenable to algebraic transformations (see \citet{Guermond2019} Sec. 7.5.2 and \Cref{ssec:idp}).

\subsubsection{Invariant Domain Preserving}\label{ssec:idp}
Invariant domain preserving schemes, introduced in \citet{Berthon2008}, are shown to preserve any convex invariant of hyperbolic systems and satisfy the discrete entropy inequality for every admissible entropy of the system \citep{Guermond2016}. These properties were achieved for a general high-order setting in \citet{Guermond2016} by using a graph-viscosity term to introduce a sufficient amount of artificial dissipation. By instead formulating the invariant domain as a set of admissible solutions, it is possible to enforce these properties via the BP-RK schemes.  

For an in-depth overview of invariant sets and domains, the reader is referred to \citet{Guermond2016} and the works therein. In brief, a set is considered to be invariant with respect to \cref{eq:gen_hype} if for any pair of states within the set, the average of the entropy solution of the Riemann problem over the Riemann fan remains within that set. From the work of \citet{Hoff1985}, it can be seen that this invariant set is convex for genuinely nonlinear hyperbolic equations. In general terms, a scheme is said to be invariant domain preserving if for some invariant set $\mathcal S$ such that $\mathbf{u}_i^{n} \in \mathcal S \ \forall \ i \in V$, the temporal update of the scheme results in $\mathbf{u}_i^{n+1} \in \mathcal S \ \forall \ i \in V$.

The formulation of invariant domain preserving bounds is shown through the example of the compressible Euler equations, written in the form of \cref{eq:gen_hype} as
\begin{equation}\label{eq:euler}
    \mathbf{u} = \begin{bmatrix}
            \rho \\ \boldsymbol{\rho v} \\ E
        \end{bmatrix}, \quad  \mathbf{F} = \begin{bmatrix}
            \boldsymbol{\rho v}\\
            \boldsymbol{\rho v}\otimes\mathbf{v} + p\mathbf{I}\\
        (E+p)\mathbf{v}
    \end{bmatrix},
\end{equation}
where $\rho$ is the density, $\boldsymbol{\rho v}$ is the momentum, $E$ is the total energy, $p = (\gamma-1)\left(E - \half\rho\|\mathbf{v}\|_2^2\right)$ is the pressure, and $\gamma > 1$ is the ratio of specific heat capacities. The symbol $\mathbf{I}$ denotes the identity matrix in $\mathbb{R}^{d\times d}$ and $\mathbf{v} = \boldsymbol{\rho v}/\rho$ denotes the velocity. From \citet{Guermond2016},  the set 
\begin{equation}
    \mathcal{S} := \{(\rho, \boldsymbol{\rho v}, E)\ |\ \rho \geq 0,\ \epsilon(\mathbf{u})  \geq 0,\ \Psi(\mathbf{u}) \geq \Psi_0\}
\end{equation}
is an invariant set for the Euler system for a specific internal energy $\epsilon(\mathbf{u})  := E/\rho - \frac{1}{2}\|\mathbf{v}\|_2^2$, specific physical entropy $\Psi(\mathbf{u})$ such that $-\Psi(\epsilon, \rho^{-1})$ is a strictly convex function, and any $\Psi_0 \in \mathbb{R}$. 

Invariant domain preserving bounds can be enforced through discrete maximum principle bounds on the density,
\begin{equation}
\rho_{\min} < \rho < \rho_{\max},
\end{equation}
and a local minimum condition on the specific physical entropy,
\begin{equation}
    \Psi(\mathbf{u}) - \Psi^{\min}_i > 0, \quad \quad \Psi(\mathbf{u}) = \frac{1}{\gamma -1} \mathrm{log} (e \rho^{-\gamma}),
\end{equation}
where $e = \rho \epsilon(\mathbf{u})$ denotes the internal energy. The minima/maxima are calculated through average Riemann solutions via the auxiliary states
\begin{equation}
    \overline{\mathbf{U}}_{ij} = \frac{1}{2} ( \mathbf{u}_i + \mathbf{u}_j) - \frac{\mathbf{c}_{ij}}{2 \lambda_{\mathrm{max}} \| \mathbf{c}_{ij} \|_2} \big ( \mathbf F (\mathbf{u}_j) - \mathbf F (\mathbf{u}_i) \big ),
\end{equation}
for $j \in \mathcal A (i)$ and some estimate of the maximum wavespeed of the system $\lambda_{\mathrm{max}}$. 

From \citet{Guermond2019}, the entropy condition can be equivalently expressed as 
\begin{equation}
    e - \rho^{\gamma} \tilde{\Psi}_i^{\min} > 0, \quad \quad \tilde{\Psi}(\mathbf{u}) := e \rho^{-\gamma} = \mathrm{exp} [(\gamma - 1) \Psi(\mathbf{u})],
\end{equation}
which yields analytic transformations for the mapping function $\mathbf G (\mathbf u )$. These constraints can be enforced through the general mappings presented in \Cref{ssec:mappings} as
\begin{align}\label{eq:eulerbounds}
    \rho &\in (\rho_{\min}, \rho_{\max}), \\
    \| \boldsymbol{\rho v} \|_2^2    & < 2\rho (E - \rho^{\gamma} \tilde{\Psi}_i^{\min}) , \\
    E &>  \rho^{\gamma} \tilde{\Psi}_i^{\min} + \frac{1}{2} \rho^{-1} \| \boldsymbol{\rho v} \|_2^2.
\end{align}
In contrast to the convex limiting approach of \citet{Guermond2019}, the mappings used to enforce these constraints have closed-form solutions and do not require the use of an iterative solver. For a one-dimensional example, the forward transformation can be given as
\begin{equation}
    \begin{bmatrix}
        w_1 
        \vphantom{\mathrm{tanh}^{-1} \left ( 2\frac{u_1-\rho_{\min}}{\rho_{\max}-\rho_{\min}} - 1 \right )}
        \\ 
        w_2 
        \vphantom{\mathrm{tanh}^{-1} \left ( u_2/\sqrt{2u_1(u_3 - \rho^{\gamma}\tilde{\Psi}_i^{\min})} \right ) }
        \\ 
        w_3 
        \vphantom{
        \log \left(u_3 - u_1\tilde{\Psi}_i^{\min} - \frac{1}{2}u_1^{-1}u_2^2 \right)}
    \end{bmatrix} =   
    G \left(
    \begin{bmatrix}
        u_1 
        \vphantom{\mathrm{tanh}^{-1} \left ( 2\frac{u_1-\rho_{\min}}{\rho_{\max}-\rho_{\min}} - 1 \right )}
        \\ 
        u_2 
        \vphantom{\mathrm{tanh}^{-1} \left ( u_2/\sqrt{2u_1(u_3 - \rho^{\gamma}\tilde{\Psi}_i^{\min})} \right ) }
        \\ 
        u_3 
        \vphantom{
        \log \left(u_3 - u_1\tilde{\Psi}_i^{\min} - \frac{1}{2}u_1^{-1}u_2^2 \right)}
    \end{bmatrix} \right)=
    \begin{bmatrix}
        \mathrm{tanh}^{-1} \left ( 2\frac{u_1-\rho_{\min}}{\rho_{\max}-\rho_{\min}} - 1 \right ) \\ 
        \mathrm{tanh}^{-1} \left ( u_2/\sqrt{2u_1(u_3 - u_1^{\gamma}\tilde{\Psi}_i^{\min})} \right ) 
        \\ 
        \log \left(u_3 - u_1^{\gamma}\tilde{\Psi}_i^{\min} - \frac{1}{2}u_2^2 / u_1 \right)
    \end{bmatrix},
\end{equation}
and the inverse transformation can be given as 
\begin{equation}
    \begin{bmatrix}
        u_1 
        \\ 
        u_2 
        \\ 
        u_3 
    \end{bmatrix} =   
    G^{-1} \left(
    \begin{bmatrix}
        w_1 
        \\ 
        w_2 
        \\ 
        w_3 
    \end{bmatrix} \right)=
    \begin{bmatrix}
        \frac{1}{2} \left(\rho_{\max}+\rho_{\min}\right) + \frac{1}{2} \left(\rho_{\min} - \rho_{\max}\right)\mathrm{tanh} \left (w_1\right ) \\ 
        \mathrm{sgn}(\zeta_2) \sqrt{2 u_1\zeta_2^2 \zeta_3/(1 - \zeta_2^2)}
        \\ 
        \zeta_3 + u_1^{\gamma}\tilde{\Psi}_i^{\min} + \frac{1}{2}u_1^{-1}u_2^2
    \end{bmatrix},
\end{equation}
where $\zeta_2 = \mathrm{tanh} (w_2)$ and $\zeta_3 = \exp(w_3)$.
Note that the reverse transform is explicit in terms of $\mathbf{w}$ as the right-hand side can be expressed in terms of just $\mathbf{w}$, although $u_1$ and $u_2$ are utilized for brevity.


\section{Implementation}\label{sec:implementation}
The BP-RK schemes were implemented and utilized on a variety of hyperbolic problems, including nonlinear problems with discontinuities. To demonstrate the potential of these schemes, a pseudospectral method was used for the spatial discretization for solutions with discontinuous features without employing any explicit shock capturing approaches. The computational domain $\Omega$ is taken to be a periodic hypercube $[0,1]^d$ (or an affine image thereof) with equispaced nodes $\{\mathbf{x}\}_{i \in V}$. The solution was approximated via a Fourier basis of degree $N-1$ as
\begin{equation}
    \mathbf{U}_h = \sum_{\| \mathbf{k} \|_\infty < N} \hat{\mathbf{u}}_k e^{2  \pi i \mathbf{k} \cdot \mathbf{x} / N},
\end{equation}
where $\mathbf{k} \in \mathbb Z^d$, $\mathbf{x} \in [0, 1]^d$, and $N \geq 1$. The projection of the nonlinear flux terms was evaluated via a pseudospectral (collocation) approach without anti-aliasing. BP-RK schemes of up to fourth order were considered in this work. These underlying RK schemes upon which the proposed methods were built are represented through the Butcher tableaux in \cref{tab:rk_schemes}.

\begin{figure}[htbp!] 
    \centering
\subfloat[RK2]{\adjustbox{width=0.12\linewidth,valign=b}{
$\renewcommand{\arraystretch}{1.5} 
\begin{array}
{c|c c }
0&  &  \\
\sfrac{1}{2} & \sfrac{1}{2} & \\
\hline
  & 0 & 1 
\end{array}
$}}
~
\subfloat[RK3]{\adjustbox{width=0.16\linewidth,valign=b}{
$\renewcommand{\arraystretch}{1.5} 
\begin{array}
{c|c c c}
0 &  &  &   \\
\sfrac{1}{3} & \sfrac{1}{3} &  &   \\
\sfrac{2}{3} & 0 & \sfrac{2}{3} &   \\
\hline
  & \sfrac{1}{4} & 0 & \sfrac{3}{4} 
\end{array}
$}}
~
\subfloat[RK4]{\adjustbox{width=0.2\linewidth,valign=b}{
$\renewcommand{\arraystretch}{1.5} 
\begin{array}
{ c|c c c c}
0 &  &  &  & \\
\sfrac{1}{2} & \sfrac{1}{2} &  &  & \\
\sfrac{1}{2} & 0 & \sfrac{1}{2} &  & \\
1 & 0 & 0 & 1 & \\
\hline
  & \sfrac{1}{6} & \sfrac{1}{3} & \sfrac{1}{3} & \sfrac{1}{6} 
\end{array}
$}}
    \caption{\label{tab:rk_schemes} Butcher tableau for Heun's RK2, Heun's RK3, and the classic RK4 method. }
\end{figure}

We consider two formulations for the bounds: discrete maximum principle (DMP) for scalar-valued solutions and invariant domain preserving (IDP) for vector-valued solutions. For both of these methods, the support $\mathcal A (i)$ at a point $\mathbf{x}_i$ is taken to be the set of direct Voronoi neighbors of $\mathbf{x}_i$, including $\mathbf{x}_i$ itself. The calculation of the maximum wavespeed for the auxiliary states in the IDP approach was performed using an exact Godunov-type Riemann solver \citep{Toro1997}.  

In comparison to standard RK schemes, the BPRK schemes can impose additional time step restrictions on the temporal integration method. These additional restrictions primarily manifest in the limit as $\mathbf{u}_i \to \partial \mathcal B_i$ as the Jacobian of the mapping tends towards infinity in this limit. Furthermore, numerical precision errors can be exacerbated for states near the boundaries as large regions of the auxiliary space may need to be mapped to a very small region of the real space. In practice, these effects were mitigated by adding a small tolerance of $\mathcal O(10^{-6})$ to the bounds, such that the maximum usable time step of the BPRK schemes was generally within 50\% of the standard RK schemes. Smaller tolerances could be used at the expense of potentially more restrictive conditions on the time step. 

To compute the mass correction step, it is necessary to evaluate the critical set distance functional $\gamma_i^* = D_{\mathcal B_i} (\mathbf{\overline{u}}^{n+1}_i, \mathbf{n})$ at each time step for all $i \in V$. We will neglect the subscript $i$ for brevity in the presentation. For scalar conservation laws with DMP bounds, this can be trivially evaluated as
\begin{equation}
    \gamma^* = 
    \begin{cases}
    u^{\max} - \overline{u}^{n+1}, \quad \mathrm {if}\ \mathrm{sgn}(\overline{{S}}) \geq 0, \\
    \overline{u}^{n+1} - u^{\min}, \quad \, \mathrm {else},
    \end{cases}
\end{equation}
utilizing that $n = \mathrm{sgn}(\overline{{S}})$ by \cref{eq:n}. For the compressible Euler equations with IDP bounds, this calculation becomes significantly more complex. In the cases where $\tilde{\Psi}^{\min} = 0$ or the specific heat ratio $\gamma$ is an integer, the set boundary can be represented as a polynomial function of the state variables, such that the distance function (i.e., the intersection of a line and a polynomial) can be computed analytically. Outside of these cases, there does not exist a closed form solution to this problem, and $\gamma^*$ must be computed numerically. 

However, let us first consider the case of $\tilde{\Psi}^{\min} = 0$, which corresponds to the set of solutions with positive density and internal energy/pressure. For the one-dimensional case, the solution is given by $\overline{\mathbf{u}}^{n+1} = [u_1, u_2, u_3] = [\rho, \rho v, E]$ and the mass defect unit vector is given by $\mathbf{n} = [n_1, n_2, n_3]$. There exists a bound for $\gamma^*$ due to the density constraint, which can be computed as 
\begin{equation}
    \gamma^{*(1)} = 
    \begin{cases}
    (\rho_{\max} - u_1)/|n_1|, \quad \mathrm {if}\ n_1 \geq 0, \\
    (u_1 - \rho_{\min})/|n_1|, \quad \, \mathrm {else}.
    \end{cases}
\end{equation}
A second bound comes from the ``entropy'' constraint (which for $\tilde{\Psi}^{\min} = 0$ is actually a constraint on the positivity of internal energy/pressure), which comes from \cref{eq:eulerbounds} by the inequality 
\begin{equation}
    u_2^2 > 2 u_1 u_3.
\end{equation}
Since the mass correction step is linear (i.e., $\mathbf{u}^{n+1} = \overline{\mathbf{u}}^{n+1} + \alpha \mathbf{n}$ for some scalar $\alpha \geq 0$), this simply becomes a quadratic constraint for which the bound can be computed as 
\begin{equation}
    \gamma^{*(2)} = \max\left [ \frac{-b \pm \sqrt{b^2 - 4ac}}{2a} \right],
\end{equation}
where
\begin{align}
    a &= n_1 n_3 - \frac{1}{2}n_2^2, \\
    b &= n_1 u_3 + u_1 n_3 - n_2 u_2, \\
    c &= u_1 u_3 - \frac{1}{2}u_1^2.
\end{align}
Note that only the positive root is considered by taking the maxima. A similar form can be given for the multi-dimensional case. The critical set distance for $\tilde{\Psi}^{\min} = 0$, which will be denoted as $\gamma^{**}$, can then be computed as
\begin{equation}
    \gamma^{**} = \min \left (\gamma^{*(1)},  \gamma^{*(2)} \right).
\end{equation}
The convexity of $\tilde{\Psi}$ and \cref{def:dist} ensure that
\begin{equation}
    0 \leq \gamma^* \leq \gamma^{**}.
\end{equation}

With this upper and lower bound on $\gamma^*$, its value can be numerically computed in an efficient manner using simple root bracketing methods such as the bisection and Illinois algorithms on a pointwise basis (i.e., independently for each $i \in V$). However, it must be noted that defining $\gamma$ utilizing $\gamma^*$ as per \cref{eq:gammai} \textit{does not form the upper bound} on the bounds preserving value of $\mathbf{S}$, and any value of $\gamma$ such that the scheme is bounds preserving and linear invariant preserving can be sufficient. In practice, it is almost always the case that utilizing $\gamma^{**}$ is sufficient, such that $\gamma$ and $\mathbf{S}$ (with the associated subscripts $i$) can be defined as 
\begin{equation}
    \gamma_i = 
      \gamma_{i}^{**} \| \overline{\mathbf{S}} \|_2/ \sum_{j \in V} m_j \gamma_{j}^{**} \quad \quad \mathrm{and} \quad \quad \mathbf{S}_i = \gamma_i \mathbf{n}.
\end{equation}
This form can be utilized to significantly improve the efficiency of the scheme as there exists an analytic solution to $\gamma$, and for the purposes of robustness, the scheme can simply check if the bounds are satisfied after the mass correction step and revert back to numerically calculating $\gamma^*$ if these bounds are violated. In this work, both methods were explored, using both the analytic approximation and a direct numerical calculation with 5 iterations of the bisection method, and the approaches were found to be virtually indistinguishable with respect to the overall results. 

In comparison to the computational cost of the spatial discretization, the cost of the mappings introduced by the BP-RK schemes was effectively negligible. Therefore, for a given time step $\Delta t$, the computational costs of the BP-RK schemes and standard RK schemes were essentially identical except for the computation of the bounds and the mass correction step. Since computing the bounds (e.g., the auxiliary states for the invariant domain) would be required for any limiting approach based on these bounds, its cost is neglected here. When the mass correction was computed analytically using only the approximate form, the resulting cost was not noticeably impacted. However, when the mass correction was computed numerically using only the exact form, the cost was significantly increased to where the evaluation of the mass correction step was on the order of the cost of the evaluation of the spatial discretization. This cost could be reduced by utilizing more efficient root finding approaches (e.g., Illinois method), lowering the iteration count, or simply utilizing the approximate form and reverting to the exact form only when necessary.

\section{Results}\label{sec:results}
\subsection{Linear Transport}
\begin{figure}[htbp!]
    \centering
    \subfloat[RK]{\label{fig:adv_error_unbounded} \adjustbox{width=0.48\linewidth,valign=b}{\begin{tikzpicture}[spy using outlines={rectangle, height=3cm,width=2.3cm, magnification=3, connect spies}]
	\begin{loglogaxis}[name=plot1,
	    axis line style={latex-latex},
		xlabel={$\Delta t$},
    	xmin=1e-5,xmax=1e-3,
    	xtick={1e-5,1e-4,1e-3},
    	ylabel={$\| e \|_2$},
    	ymin=1e-16,ymax=1e-1,
    	grid=both,
    	x dir=reverse,
    	legend style={at={(0.97,0.03)},anchor=south east,font=\small},
    	legend cell align={left},
    	style={font=\normalsize}]
    
		\addplot[color={black}, style={thick},  forget plot]
		coordinates{(5e-5,0.3*0.000876677) 
		(2.5e-5,0.3*0.000876677/2)
		(5e-5,0.3*0.000876677/2)
		(5e-5,0.3*0.000876677)};
    	\addplot[color={black}, style={very thick}, mark=triangle, mark options={solid, scale=1.3}] coordinates {
    	    (1.00E-04, 0.001762049)
    		(5.00E-05, 0.000876677)
    		(2.50E-05, 0.000437257)
    		(1.25E-05, 0.000218359)
    	};
		\addlegendentry{Euler};
		\node at (axis cs:3.75e-5,0.00001972523) {1};
		
		\addplot[color={black}, style={thick},  forget plot]
		coordinates{(1e-4,0.3*3.65E-07) 
		(5e-5,0.3*3.65E-07/4)
		(1e-4,0.3*3.65E-07/4)
		(1e-4,0.3*3.65E-07)};
    	\addplot[color={black}, style={very thick, dashed}, mark=square, mark options={solid, scale=1.0}] coordinates {
    		(2.00E-04, 1.46E-06)
    		(1.00E-04, 3.65E-07)
    		(5.00E-05, 9.14E-08)
    		(2.50E-05, 2.28E-08)
    	};
		\addlegendentry{RK2};
		\node at (axis cs:7.5e-5,4.10625e-9) {2};
		
		\addplot[color={black}, style={thick},  forget plot]
		coordinates{(2e-4,0.3*4.59E-10) 
		(1e-4,0.3*4.59E-10/8)
		(2e-4,0.3*4.59E-10/8)
		(2e-4,0.3*4.59E-10)};
    	\addplot[color={black}, style={very thick, dotted}, mark=diamond, mark options={solid, scale=1.4}] coordinates {
    		(4.00E-04, 3.67E-09)
    		(2.00E-04, 4.59E-10)
    		(1.00E-04, 5.74E-11)
    		(5.00E-05, 7.18E-12)
    	};
		\addlegendentry{RK3};
		\node at (axis cs:1.5e-4,2.581875e-12) {3};
		
		\addplot[color={black}, style={thick},  forget plot]
		coordinates{(4e-4,0.3*1.85E-12) 
		(2e-4,0.3*1.85E-12/16)
		(4e-4,0.3*1.85E-12/16)
		(4e-4,0.3*1.85E-12)};
    	\addplot[color={black}, style={very thick, dash dot}, mark=o, mark options={solid, scale=1.0}] coordinates {
    		(8.00E-04, 2.95E-11)
    		(4.00E-04, 1.85E-12)
    		(2.00E-04, 1.14E-13)
    		(1.00E-04, 1.54E-14)
    	};
		\addlegendentry{RK4};
		\node at (axis cs:3e-4,5.203125e-15) {4};
		

	\end{loglogaxis}
\end{tikzpicture}}}
    ~
    \subfloat[BP-RK]{\label{fig:adv_error_bounded} \adjustbox{width=0.48\linewidth,valign=b}{\begin{tikzpicture}[spy using outlines={rectangle, height=3cm,width=2.3cm, magnification=3, connect spies}]
	\begin{loglogaxis}[name=plot1,
		axis line style={latex-latex},
		xlabel={$\Delta t$},
    	xmin=1e-5,xmax=1e-3,
    	xtick={1e-5,1e-4,1e-3},
    	ylabel={$\| e \|_2$},
    	ymin=1e-16,ymax=1e-1,
    	grid=both,
    	x dir=reverse,
    	legend style={at={(0.97,0.03)},anchor=south east,font=\small},
    	legend cell align={left},
    	style={font=\normalsize}]
    	
		\addplot[color={black}, style={thick}, forget plot
        ]
		coordinates{(5e-5,0.3*0.000798814) 
		(2.5e-5,0.3*0.000798814/2)
		(5e-5,0.3*0.000798814/2)
		(5e-5,0.3*0.000798814)};
    	\addplot[color={black}, style={very thick}, mark=triangle, mark options={solid, scale=1.3}] coordinates {
    	    (1.00E-04, 0.001606245)
    		(5.00E-05, 0.000798814)
    		(2.50E-05, 0.000398545)
    		(1.25E-05, 0.000199059)
    	};
		\node at (axis cs:3.75e-5,0.00001797331 ) {1};

		\addplot[color={black}, style={thick}, forget plot
        ]
		coordinates{(1e-4,0.3*3.66E-06) 
		(5e-5,0.3*3.66E-06/4)
		(1e-4,0.3*3.66E-06/4)
		(1e-4,0.3*3.66E-06)};
    	\addplot[color={black}, style={very thick, dashed}, mark=square, mark options={solid, scale=1.0}] coordinates {
    		(2.00E-04, 1.46E-05)
    		(1.00E-04, 3.66E-06)
    		(5.00E-05, 9.14E-07)
    		(2.50E-05, 2.29E-07)
    	};
		\node at (axis cs:7.5e-5,4.1175e-8) {2};
		
		\addplot[color={black}, style={thick}, forget plot
        ]
		coordinates{(2e-4,0.3*1.10E-08) 
		(1e-4,0.3*1.10E-08/8)
		(2e-4,0.3*1.10E-08/8)
		(2e-4,0.3*1.10E-08)};
    	\addplot[color={black}, style={very thick, dotted}, mark=diamond, mark options={solid, scale=1.4}] coordinates {
    		(4.00E-04, 8.80E-08)
    		(2.00E-04, 1.10E-08)
    		(1.00E-04, 1.37E-09)
    		(5.00E-05, 1.71E-10)
    	};
		\node at (axis cs:1.5e-4,6.1875e-11 ) {3};

		\addplot[color={black}, style={thick}, forget plot
        ]
		coordinates{(4e-4,0.3*6.81E-10) 
		(2e-4,0.3*6.81E-10/16)
		(4e-4,0.3*6.81E-10/16)
		(4e-4,0.3*6.81E-10)};
    	\addplot[color={black}, style={very thick, dash dot}, mark=o, mark options={solid, scale=1.0}] coordinates {
    		(8.00E-04, 1.08E-08)
    		(4.00E-04, 6.81E-10)
    		(2.00E-04, 4.75E-11)
    		(1.00E-04, 3.08E-12)
    	};
		\node at (axis cs:3e-4,1.9153125e-12) {4};

	\end{loglogaxis}
\end{tikzpicture}}}
    \caption{\label{fig:adv_error} Convergence of the $L^2$ norm of the error with respect to time step using Runge-Kutta (left) and bounds preserving Runge-Kutta temporal integration (right). Solid triangles represent the theoretical convergence rates.}
\end{figure}
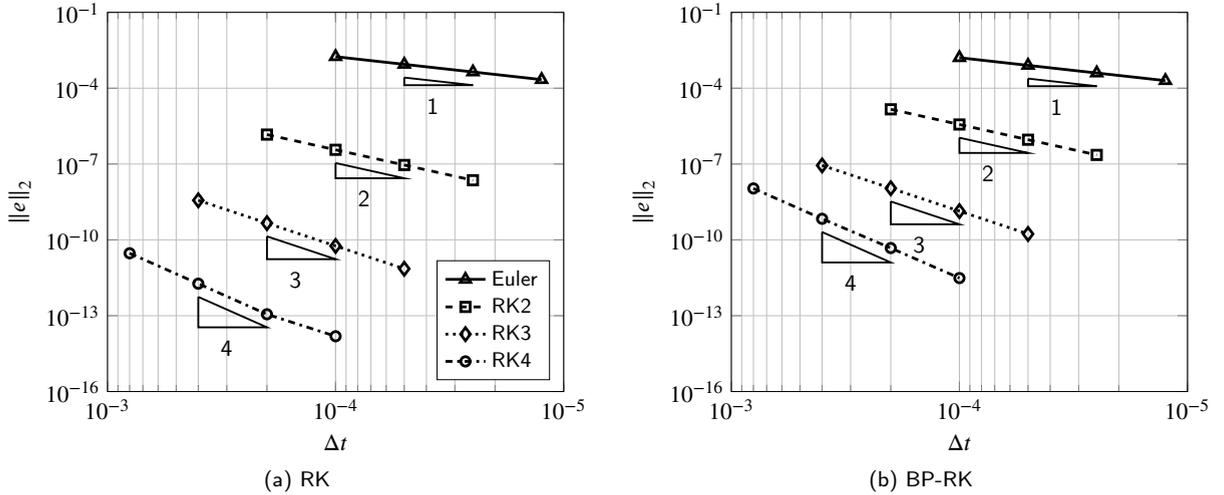

The convergence rates and bounds-preserving properties of the BP-RK schemes were tested on the linear transport equation in one dimension:
\begin{equation}
    \partial_t u + \partial_x (cu) = 0.
\end{equation}
The transport velocity was set as $c = 1$ and the domain was set as $\Omega = [0,1]$ with periodic boundary conditions. A smooth sinusoidal initial condition, given as 
\begin{equation}
    u_0(x) = \mathrm{sin} (2 \pi x),
\end{equation}
was used to validate the results of \cref{thm:convergence}. The $L^2$ norm of the solution error, defined as 
\begin{equation}
    \| e \|_2 = \sqrt{\sum_{i=0}^{N-1} h \big( u_i - u_0(x_i) \big )^2},
\end{equation}
with $h = 1/N$, was computed at $t = 100$, corresponding to 100 traverses through the domain. The convergence rate of the $L^2$ error with respect to $\Delta t$ for a pseudospectral spatial discretization with $N = 32$ using RK and BP-RK temporal integration with DMP bounds is shown in \cref{fig:adv_error} and tabulated in \cref{tab:l2_error}. The BP-RK schemes recovered the theoretical convergence rates of the underlying RK schemes, as stated by \cref{thm:convergence}, albeit with larger leading error constants for $p > 1$. These effects may be attributed to the larger impact of numerical precision errors due to the nonlinear transformations required by the BP-RK schemes. 

\begin{table}[htbp!]
    \centering
    \begin{tabular}{ccc}
        \toprule
        Order & RK & BP-RK\\
        \midrule
        1 & 1.004 & 1.005\\
        2 & 2.000 & 2.000\\
        3 & 2.999 & 3.002\\
        4 & 3.674 & 3.959\\\bottomrule
    \end{tabular}
    \caption{Convergence rates of the $L^2$ norm of the error with respect to time step using Runge-Kutta and bounds preserving Runge-Kutta temporal integration.}
    \label{tab:l2_error}
\end{table}

The ability of the BP-RK schemes to enforce bounds was initially evaluated for a linear transport problem with a non-smooth initial condition given by 
\begin{equation}
    u_0(x) = \begin{cases}
        \exp\left( -300(2x - 0.3)^2 \right), &\mbox{if } |2x - 0.3| \leq 0.25, \\
        1, &\mbox{if } |2x - 0.9| \leq 0.2, \\
        \sqrt{1 - \left(\frac{2x-1.6}{0.2}\right)^2}, &\mbox{if } |2x - 1.6| \leq 0.2, \\
        0, &\mbox{else}.
    \end{cases} 
\end{equation}
\begin{figure}[htbp!]
    \centering
    \subfloat[RK4]{\label{fig:dadv_unbounded} \adjustbox{width=0.48\linewidth,valign=b}{    \begin{tikzpicture}[spy using outlines={rectangle, height=3cm,width=2.3cm, magnification=3, connect spies}]
		\begin{axis}[name=plot1,
		    axis line style={latex-latex},
		    axis x line=left,
            axis y line=left,
            clip mode=individual,
		    xlabel={$x$},
		    xtick={0,0.2,0.4,0.6,0.8,1},
    		xmin=0,
    		xmax=1,
    		x tick label style={
        		/pgf/number format/.cd,
            	fixed,
            	fixed zerofill,
            	precision=1,
        	    /tikz/.cd},
    		ylabel={$u$},
    		ylabel style={rotate=-90},
    		ytick={-0.5, 0, 0.5, 1.0, 1.5},
    		ymin=-0.5,
    		ymax=1.5,
    		y tick label style={
        		/pgf/number format/.cd,
            	fixed,
            	fixed zerofill,
            	precision=1,
        	    /tikz/.cd},
    		legend style={at={(0.5,1)},anchor=north ,font=\small, column sep=0.2cm},
    		legend cell align={left},
    		legend columns=3,
    		style={font=\normalsize}]
    		
			\addplot[color=gray, style={ultra thin}, only marks, mark=o, mark options={scale=0.8}, mark repeat = 4, mark phase = 0]
				table[x=x,y=u,col sep=comma,unbounded coords=jump]{./figs/data/shapes_exact.csv};
    		\addlegendentry{Exact}
    		
			\addplot[color={black}, style={very thick}]
				table[x=x,y=u,col sep=comma,unbounded coords=jump]{./figs/data/advection_shapes_n128_unbounded.csv};
    		\addlegendentry{$N = 128$}
    		
			\addplot[color={red!90!black}, style={ultra thick, dotted}]
				table[x=x,y=u,col sep=comma,unbounded coords=jump]{./figs/data/advection_shapes_n256_unbounded.csv};
    		\addlegendentry{$N = 256$}

		\end{axis}

	\end{tikzpicture}}}
    ~
    \subfloat[BP-RK4]{\label{fig:dadv_bounded} \adjustbox{width=0.48\linewidth,valign=b}{    \begin{tikzpicture}[spy using outlines={rectangle, height=3cm,width=2.3cm, magnification=3, connect spies}]
		\begin{axis}[name=plot1,
		    axis x line=left,
            axis y line=left,
		    xlabel={$x$},
		    xtick={0,0.2,0.4,0.6,0.8,1},
    		xmin=0,
    		xmax=1,
    		x tick label style={
        		/pgf/number format/.cd,
            	fixed,
            	fixed zerofill,
            	precision=1,
        	    /tikz/.cd},
    		ylabel={$u$},
    		ylabel style={rotate=-90},
    		ytick={-0.5, 0, 0.5, 1.0, 1.5},
    		ymin=-0.5,
    		ymax=1.5,
    		y tick label style={
        		/pgf/number format/.cd,
            	fixed,
            	fixed zerofill,
            	precision=1,
        	    /tikz/.cd},
    		legend style={at={(0.5,1)},anchor=north ,font=\small, column sep=0.2cm},
    		legend cell align={left},
    		legend columns=3,
    		style={font=\normalsize}]
    		
			\addplot[color=gray, style={ultra thin}, only marks, mark=o, mark options={scale=0.8}, mark repeat = 4, mark phase = 0]
				table[x=x,y=u,col sep=comma,unbounded coords=jump]{./figs/data/shapes_exact.csv};

			\addplot[color={black}, style={very thick}]
				table[x=x,y=u,col sep=comma,unbounded coords=jump]{./figs/data/advection_shapes_n128_bounded.csv};
    		
			\addplot[color={red!90!black}, style={ultra thick, dotted}]
				table[x=x,y=u,col sep=comma,unbounded coords=jump]{./figs/data/advection_shapes_n256_bounded.csv};

		\end{axis}

	\end{tikzpicture}}}
    \caption{\label{fig:dadv} Solution of the linear transport equation at $t = 100$ using RK4 (left) and bounds preserving RK4 temporal integration (right), pseudospectral spatial discretization ($N = 128$, $256$), and discrete maximum principle bounds. }
\end{figure}
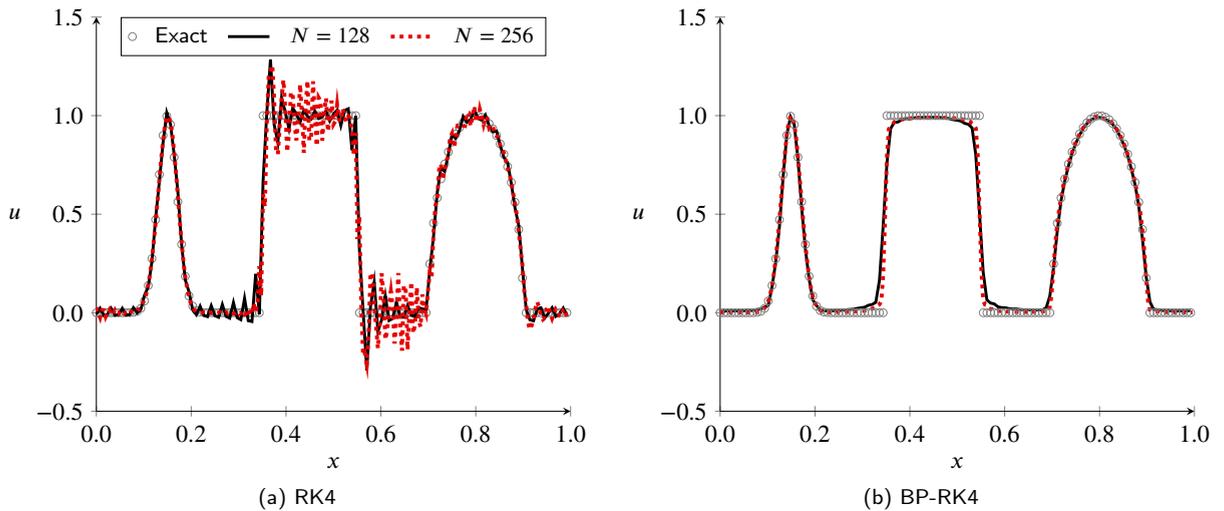
The solution at $t=100$ using a pseudospectral discretization with $N = 128$, $256$ is shown in \cref{fig:dadv} for RK4 and BP-RK4 temporal integration with DMP bounds with $\Delta t = 4{\cdot}10^{-3}$ and $2{\cdot}10^{-3}$, respectively. Without the BP-RK schemes, the solution was highly oscillatory with significant overshoots and undershoots around discontinuous features. When the BP-RK4 was used, the predicted solution was in excellent agreement with the exact solution, and no spurious oscillations were observed. Due to the enforcement of the discrete maximum principle, the solution remained within the range of the initial conditions. 

\subsection{Nonlinear Transport}
The BP-RK schemes were then evaluated for nonlinear hyperbolic conservation laws that develop discontinuities from smooth initial conditions. The inviscid Burgers equation in one-dimension, given by 
\begin{equation}
    \partial_t u + \partial_x \big(\frac{1}{2}u^2 \big) = 0,
\end{equation}
with the initial conditions
\begin{equation}
    u_0(x) = \mathrm{sin} (2 \pi x) + 2,
\end{equation}
was solved on the periodic domain $\Omega = [0,1]$ using a pseudospectral spatial discretization with $N = 64$. 
\begin{figure}[htbp!]
    \centering 
    \subfloat[Solution]{
    \adjustbox{width=0.45\linewidth,valign=b}{    \begin{tikzpicture}[spy using outlines={rectangle, height=3cm,width=2.3cm, magnification=3, connect spies}]
		\begin{axis}[name=plot1,
		    axis line style={latex-latex},
		    axis x line=left,
            axis y line=left,
            clip mode=individual,
		    xlabel={$x$},
		    xtick={0,0.2,0.4,0.6,0.8,1},
    		xmin=0,
    		xmax=1,
    		x tick label style={
        		/pgf/number format/.cd,
            	fixed,
            	fixed zerofill,
            	precision=1,
        	    /tikz/.cd},
    		ylabel={$u$},
    		ylabel style={rotate=-90},
    		ytick={0.5, 1.5, 2.5, 3.5},
    		ymin=0.5,
    		ymax=3.5,
    		y tick label style={
        		/pgf/number format/.cd,
            	fixed,
            	fixed zerofill,
            	precision=1,
        	    /tikz/.cd},
    		legend style={at={(1,1)},anchor=north east,font=\small},
    		legend cell align={left},
    		style={font=\normalsize}]
    		
			\addplot[color=gray, style={ultra thin}, only marks, mark=o, mark options={scale=0.8}, mark repeat = 1, mark phase = 0]
				table[x=x,y=u,col sep=comma,unbounded coords=jump]{./figs/data/sine_exact.csv};
				
    		\addlegendentry{$u_0$}
    		

			\addplot[color={black}, style={very thick}]
				table[x=x,y=u,col sep=comma,unbounded coords=jump]{./figs/data/burgers_n64_bounded.csv};
    		\addlegendentry{BP-RK2}
    		
			\addplot[color={red!90!black}, style={ultra thick, dotted}]
				table[x=x,y=u,col sep=comma,unbounded coords=jump]{./figs/data/burgers_n64_unbounded.csv};
    		\addlegendentry{RK2}

		\end{axis}

	\end{tikzpicture}}}
    \subfloat[Mass Residual]{
    \adjustbox{width=0.51\linewidth,valign=b}{    \begin{tikzpicture}[spy using outlines={rectangle, height=3cm,width=2.3cm, magnification=3, connect spies}]
		\begin{axis}
		    [name=plot1,
		    axis line style={latex-latex},
		    axis x line=left,
            axis y line=left,
            clip mode=individual,
		    xlabel={$t$},
		    xtick={0,0.2,0.4,0.6,0.8,1},
    		xmin=0,
    		xmax=1,
    		x tick label style={
        		/pgf/number format/.cd,
            	fixed,
            	fixed zerofill,
            	precision=1,
        	    /tikz/.cd},
    		ylabel={$\langle u(t) \rangle - \langle u_0 \rangle$},
    		ytick={-2e-15, -1e-15, 0, 1e-15, 2e-15},
    		ymin=-2e-15,
    		ymax=2e-15,
    		y tick label style={
        		/pgf/number format/.cd,
            	fixed,
            	fixed zerofill,
            	precision=1,
        	    /tikz/.cd},
    		legend style={at={(0.03,1)},anchor=north west,font=\small},
    		legend cell align={left},
    		style={font=\normalsize}]
    		
			\addplot[color={black}, style={very thick}, each nth point={10}]
				table[x index=0, y index=1,col sep=comma,unbounded coords=jump]{./figs/data/burgers_res.csv};
    		\addlegendentry{Mass residual}       
    		
    		\addlegendimage{color={black}, style={very thick}, dashed}
    		\addlegendentry{$\overline{S}$ magnitude}       

		\end{axis} 		
			
		\begin{axis}[name=plot2,
		    axis line style={latex-latex},
            axis y line=right,
            axis x line=none,
            clip mode=individual,
		    xtick={0,0.2,0.4,0.6,0.8,1},
    		xmin=0,
    		xmax=1,
    		x tick label style={
        		/pgf/number format/.cd,
            	fixed,
            	fixed zerofill,
            	precision=1,
        	    /tikz/.cd},
    		ylabel={$|\overline{S}|$},
    		ylabel style={rotate=-90},
    		ytick={0, 2e-5, 4e-5, 6e-5, 8e-5, 1e-4},
    		ymin=0,
    		ymax=1e-4,
    		y tick label style={
        		/pgf/number format/.cd,
            	fixed,
            	fixed zerofill,
            	precision=1,
        	    /tikz/.cd},
    		style={font=\normalsize}]
    		
			\addplot[color={black}, style={very thick}, dashed]
				table[x index=0,y index=1,col sep=comma,unbounded coords=jump]{./figs/data/burgers_s.csv};
		\end{axis}

	\end{tikzpicture}}}
    \caption{\label{fig:burgers} \textit{Left}: Solution of the Burgers equation at $t = 1$ using RK2 and bounds preserving RK2 temporal integration, pseudospectral spatial discretization ($N = 64$), and discrete maximum principle bounds. \textit{Right}: Space-integrated mass residual (solid) and magnitude of $\overline{S}$ (dashed) over time. }
\end{figure}
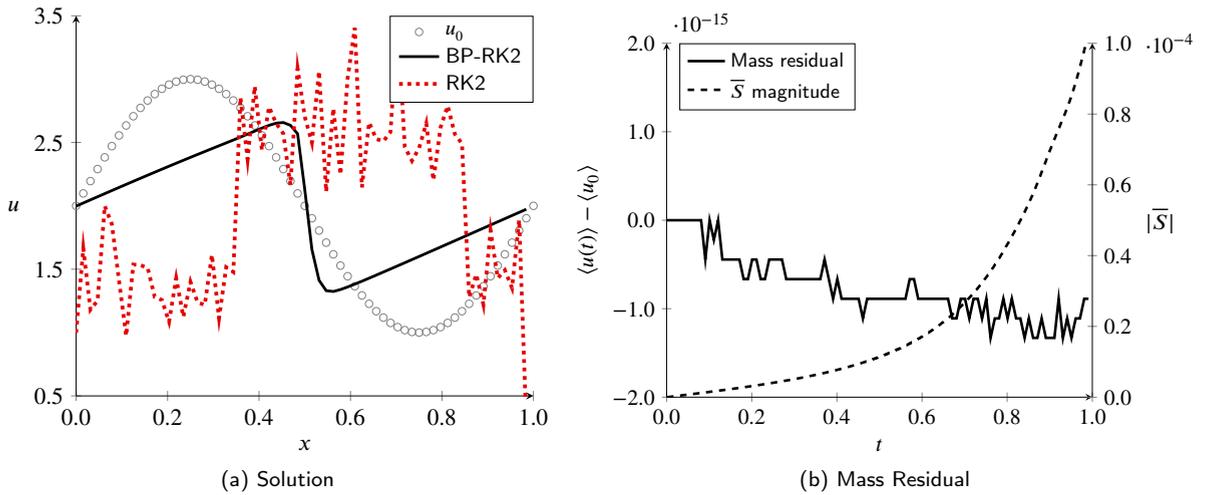
The computed solution at $t = 1$ is shown in \cref{fig:burgers} using RK2 and BP-RK2 temporal integration with DMP bounds with $\Delta t = 1{\cdot}10^{-3}$. Without the BP-RK schemes, the presence of Gibbs phenomena compounded with the nonlinearities in the transport equation made the results unusable for practical purposes due to their highly oscillatory nature. With the BP-RK2 scheme and DMP bounds, the solution remained well-behaved even as discontinuities became present in the solution through the introduction of an adequate amount of numerical dissipation around the shock via the temporal scheme. \Cref{fig:burgers} additionally shows the space-integrated mass of the solution, $\langle u(t) \rangle = \sum_{i \in V} m_i u_i(t)$, and the magnitude of the $\overline{S}$ term with respect to time. The scheme conserved mass down to numerical precision, as expected by \cref{thm:lin_inv}. The effect of the BP-RK2 scheme in enforcing the bounds is represented through the magnitude of the $\overline{S}$ term, as larger values of $\overline{S}$ indicate larger deviations from the underlying RK scheme. The magnitude of $\overline{S}$ was initially low as the solution was smooth and therefore the mappings were approximately linear, but as the solution began to develop a discontinuity, the magnitude increased to compensate for the mass defect due to the nonlinearities in the mapping functions. However, even at its maximum, this defect was still orders of magnitude lower than the overall mass of the system.

\subsection{Euler Equations}
The assessment of the BP-RK schemes was extended to vector-valued solutions and to higher-dimensions through the Euler equations, as presented in \cref{eq:euler}. For brevity, we express the solution in terms of the primitive variables $\mathbf{q}=[\rho,\mathbf{v},P]^T$.

\subsubsection{Sod Shock Tube}
The Sod shock tube problem \citep{Sod1978} was used to evaluate the ability of the proposed scheme to predict the three main features of Riemann problems: expansion fans, contact discontinuities, and shock waves. The problem is defined on the domain $\Omega =[0,1]$ with the initial conditions
\begin{equation}
    \mathbf{q}_0(x) = \begin{cases}
        \mathbf{q}_l, &\mbox{if } x \leq 0.5,\\
        \mathbf{q}_r, &\mbox{else},
    \end{cases} \quad \mathrm{given} \quad \mathbf{q}_l = \begin{bmatrix}
        1 \\ 0 \\ 1
    \end{bmatrix}, \quad \mathbf{q}_r = \begin{bmatrix}
        0.125 \\ 0 \\ 0.1
    \end{bmatrix}.
\end{equation}
We further consider the modified form of the Sod shock tube presented as Test 1 in \citet{Toro1997}, given by the initial states
\begin{equation}
\mathbf{q}_l = \begin{bmatrix}
        1 \\ 0.75 \\ 1
    \end{bmatrix}, \quad \mathbf{q}_r = \begin{bmatrix}
        0.125 \\ 0 \\ 0.1
    \end{bmatrix}.
\end{equation}
This form exhibits a sonic point in the rarefaction wave and is useful for testing the entropy satisfaction ability of numerical schemes.

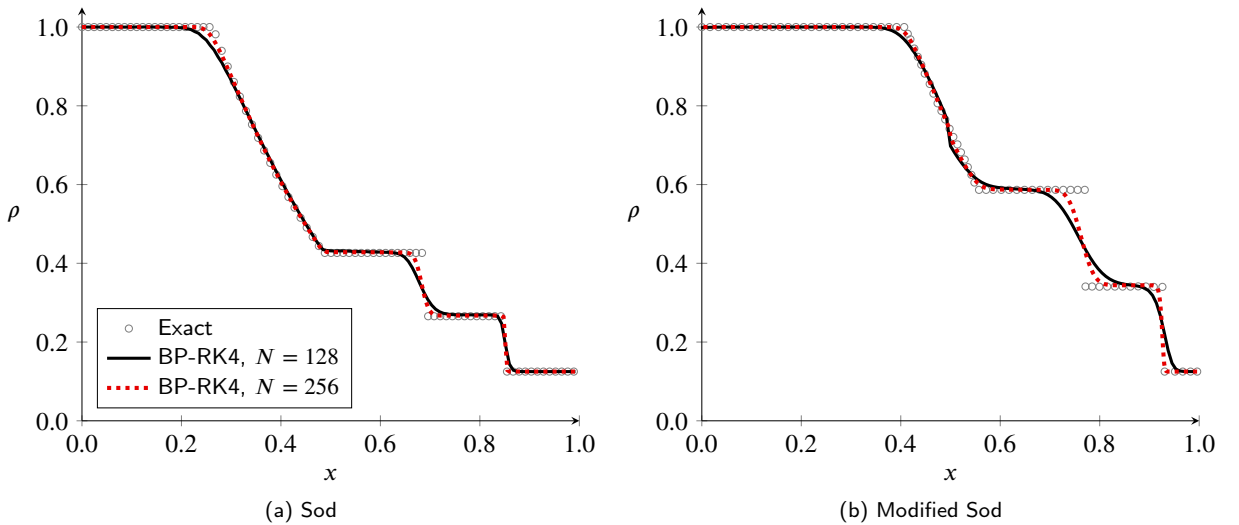
\begin{figure}[htbp!]
    \centering
    
    \subfloat[Sod]{
    \adjustbox{width=0.49\linewidth,valign=b}{    \begin{tikzpicture}[spy using outlines={rectangle, height=3cm,width=2.3cm, magnification=3, connect spies}]
		\begin{axis}[name=plot1,
		    axis line style={latex-latex},
		    axis x line=left,
            axis y line=left,
            clip mode=individual,
		    xlabel={$x$},
		    xtick={0,0.2,0.4,0.6,0.8,1},
    		xmin=0,
    		xmax=1,
    		x tick label style={
        		/pgf/number format/.cd,
            	fixed,
            	fixed zerofill,
            	precision=1,
        	    /tikz/.cd},
    		ylabel={$\rho$},
    		ylabel style={rotate=-90},
    		ytick={0,0.2,0.4,0.6,0.8,1},
    		ymin=0,
    		ymax=1.05,
    		y tick label style={
        		/pgf/number format/.cd,
            	fixed,
            	fixed zerofill,
            	precision=1,
        	    /tikz/.cd},
    		legend style={at={(0.03,0.03)},anchor=south west,font=\small},
    		legend cell align={left},
    		style={font=\normalsize}]
    		
			\addplot[color=gray, style={thin}, only marks, mark=o, mark options={scale=0.7}, mark repeat = 5, mark phase = 0]
				table[x=x,y=r,col sep=comma,unbounded coords=jump]{./figs/data/sod_exact_reference.csv};
    		\addlegendentry{Exact}
    		
			\addplot[color={black}, style={very thick}]
				table[x=x,y=r,col sep=comma,unbounded coords=jump]{./figs/data/sod_n128_exact.csv};
    		\addlegendentry{BP-RK4, $N = 128$}
    		
			\addplot[color={red!90!black}, style={ultra thick, dotted}]
				table[x=x,y=r,col sep=comma,unbounded coords=jump]{./figs/data/sod_n256_exact.csv};
    		\addlegendentry{BP-RK4, $N = 256$}

		\end{axis}

	\end{tikzpicture}}}
    \subfloat[Modified Sod]{
    \adjustbox{width=0.49\linewidth,valign=b}{    \begin{tikzpicture}[spy using outlines={rectangle, height=3cm,width=2.3cm, magnification=3, connect spies}]
		\begin{axis}[name=plot1,
		    axis line style={latex-latex},
		    axis x line=left,
            axis y line=left,
            clip mode=individual,
		    xlabel={$x$},
		    xtick={0,0.2,0.4,0.6,0.8,1},
    		xmin=0,
    		xmax=1,
    		x tick label style={
        		/pgf/number format/.cd,
            	fixed,
            	fixed zerofill,
            	precision=1,
        	    /tikz/.cd},
    		ylabel={$\rho$},
    		ylabel style={rotate=-90},
    		ytick={0,0.2,0.4,0.6,0.8,1},
    		ymin=0,
    		ymax=1.05,
    		y tick label style={
        		/pgf/number format/.cd,
            	fixed,
            	fixed zerofill,
            	precision=1,
        	    /tikz/.cd},
    		legend style={at={(0.03,0.03)},anchor=south west,font=\small},
    		legend cell align={left},
    		style={font=\normalsize}]
    		
			\addplot[color=gray, style={thin}, only marks, mark=o, mark options={scale=0.7}, mark repeat = 1, mark phase = 0]
				table[x=x,y=rho,col sep=comma,unbounded coords=jump]{./figs/data/modsod_ref.csv};

			\addplot[color={black}, style={very thick}]
				table[x=x,y=rho,col sep=comma,unbounded coords=jump]{./figs/data/modsod_n1.csv};
    		
			\addplot[color={red!90!black}, style={ultra thick, dotted}]
				table[x=x,y=rho,col sep=comma,unbounded coords=jump]{./figs/data/modsod_n2.csv};

		\end{axis}

	\end{tikzpicture}}}
    \caption{\label{fig:sod} Density profile of the Sod shock tube (left) and modified Sod shock tube (right) at $t = 0.2$ using bounds preserving RK4 temporal integration, pseudospectral spatial discretization ($N = 128$, $256$), and invariant domain preserving bounds. }
\end{figure}

The problem was periodized by extending the domain to $\Omega =[-1,1]$ and reflecting the initial conditions about $x = 0$. The resolution of the spatial scheme $N$ is presented in terms of the half-domain. The density profile of both solutions at $t = 0.2$ using BP-RK4 temporal integration, pseudospectral spatial discretization ($N = 128$, $256$), and invariant domain preserving bounds with $\Delta t = 1{\cdot}10^{-3}$ and $5{\cdot}10^{-4}$, respectively, are shown in \cref{fig:sod}. Given the coarse resolution ($N = 128$), the results were in good agreement with the exact solutions, showing no observable spurious oscillations and good resolution of the contact discontinuity, shock wave, and expansion fan without any explicit shock capturing approaches. When the resolution was increased ($N = 256$), even better agreement was observed, with particularly notable improvements at the front of the expansion fan and the contact discontinuity. Without the BP-RK scheme, the solutions diverged due to the spurious oscillations causing negative density and pressure values. 

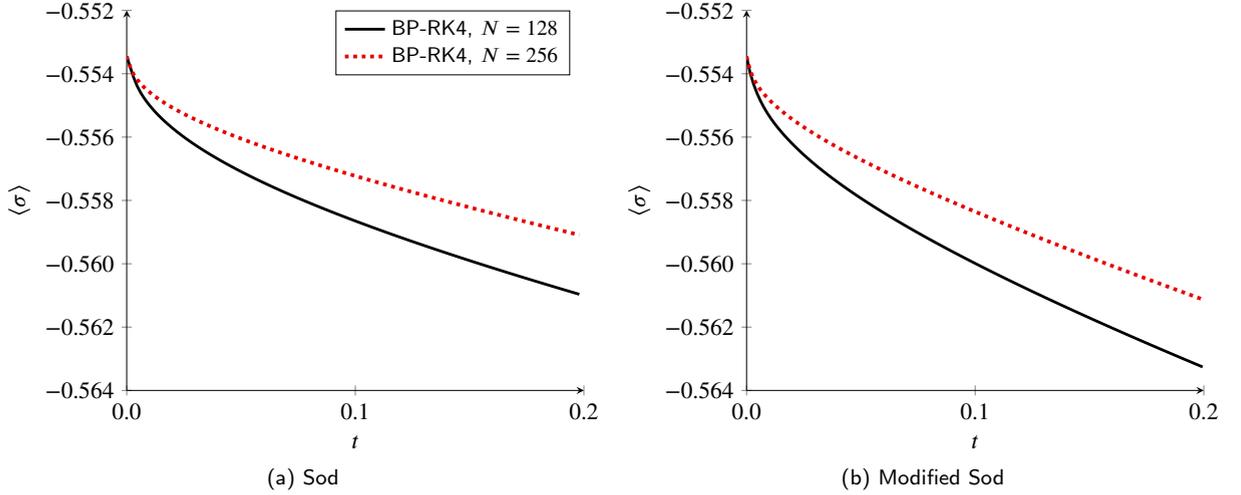
\begin{figure}[htbp!]
    \centering
    \subfloat[Sod]{
    \adjustbox{width=0.49\linewidth,valign=b}{    \begin{tikzpicture}[spy using outlines={rectangle, height=3cm,width=2.3cm, magnification=3, connect spies}]
		\begin{axis}[name=plot1,
		    axis line style={latex-latex},
		    axis x line=left,
            axis y line=left,
            clip mode=individual,
		    xlabel={$t$},
		    xtick={0,0.1, 0.2},
    		xmin=0,
    		xmax=0.2,
    		x tick label style={
        		/pgf/number format/.cd,
            	fixed,
            	fixed zerofill,
            	precision=1,
        	    /tikz/.cd},
    		ylabel={$\langle \sigma \rangle$},
    		ymin=-0.564,
    		ymax=-0.552,
    		y tick label style={
        		/pgf/number format/.cd,
            	fixed zerofill,
            	precision=3,
        	    /tikz/.cd        	    },
    		legend style={at={(0.97,1)},anchor=north east,font=\small},
    		legend cell align={left},
    		style={font=\normalsize}]
    		
			\addplot[color={black}, style={very thick}]
				table[x index=0,y index=1,col sep=comma,unbounded coords=jump]{./figs/data/dent_sod128.csv};
    		\addlegendentry{BP-RK4, $N = 128$}
    
			\addplot[color={red!90!black}, style={ultra thick, dotted}]
				table[x index=0,y index=1,col sep=comma,unbounded coords=jump]{./figs/data/dent_sod256.csv};
    		\addlegendentry{BP-RK4, $N = 256$}
		\end{axis}
	\end{tikzpicture}}}
    \subfloat[Modified Sod]{
    \adjustbox{width=0.49\linewidth,valign=b}{    \begin{tikzpicture}[spy using outlines={rectangle, height=3cm,width=2.3cm, magnification=3, connect spies}]
		\begin{axis}[name=plot1,
		    axis line style={latex-latex},
		    axis x line=left,
            axis y line=left,
            clip mode=individual,
		    xlabel={$t$},
		    xtick={0,0.1, 0.2},
    		xmin=0,
    		xmax=0.2,
    		x tick label style={
        		/pgf/number format/.cd,
            	fixed,
            	fixed zerofill,
            	precision=1,
        	    /tikz/.cd},
    		ylabel={$\langle \sigma \rangle$},
    		ymin=-0.564,
    		ymax=-0.552,
    		y tick label style={
        		/pgf/number format/.cd,
            	fixed zerofill,
            	precision=3,
        	    /tikz/.cd        	    },
    		legend style={at={(0.97,1)},anchor=north east,font=\small},
    		legend cell align={left},
    		style={font=\normalsize}]
    		
			\addplot[color={black}, style={very thick}]
				table[x index=0,y index=1,col sep=comma,unbounded coords=jump]{./figs/data/dent_modsod128.csv};
    
			\addplot[color={red!90!black}, style={ultra thick, dotted}]
				table[x index=0,y index=1,col sep=comma,unbounded coords=jump]{./figs/data/dent_modsod256.csv};
		\end{axis}
	\end{tikzpicture}}}
    \caption{\label{fig:sod_dent} Domain integrated numerical entropy for the Sod shock tube (left) and modified Sod shock tube (right) using bounds preserving RK4 temporal integration, pseudospectral spatial discretization ($N = 128$, $256$), and invariant domain preserving bounds.}
\end{figure}

For the modified problem, the \textit{entropy glitch} \citep{Toro1997} is present at the sonic point, shown as a discontinuity in the density profile in the rarefaction wave. This behavior is expected as the bounds were computed using a Godunov-type approach, and the entropy glitch is commonplace in Godunov methods (among others) even though they are theoretically entropy satisfying. However, as the resolution was increased, the size of the discontinuity decreased which is expected if the approach is entropy satisfying. 
To verify the entropy satisfying properties of the approach, the temporal evolution of the domain integrated numerical entropy $\sigma = -\rho \log (\Psi)$ is shown in \cref{fig:sod_dent} for both the Sod and modified Sod cases. It can be seen that for both cases and resolutions, the methods are strictly entropy dissipative, with a slight dependency on the spatial resolution due to the increased numerical dissipation introduced by the invariant domain preserving bounds at lower resolutions.

\begin{figure}[htbp!]
    \centering
    \subfloat[Sod]{
    \adjustbox{width=0.49\linewidth,valign=b}{    \begin{tikzpicture}[spy using outlines={rectangle, height=3cm,width=2.3cm, magnification=3, connect spies}]
		\begin{axis}[name=plot1,
		    axis line style={latex-latex},
		    axis x line=left,
            axis y line=left,
            clip mode=individual,
		    xlabel={$x$},
		    xtick={0,0.2,0.4,0.6,0.8,1},
    		xmin=0,
    		xmax=1,
    		x tick label style={
        		/pgf/number format/.cd,
            	fixed,
            	fixed zerofill,
            	precision=1,
        	    /tikz/.cd},
    		ylabel={$\gamma^*$, \textcolor{red!90!black}{$\gamma^{**}$}},
    		ytick={0.00, 0.01, 0.02, 0.03, 0.04, 0.05},
    		ymax=0.05,
    		y tick label style={
        		/pgf/number format/.cd,
        	    /tikz/.cd        	    },
    		legend style={at={(0.03,1)},anchor=north west,font=\small},
    		legend cell align={left},
    		style={font=\normalsize}]
    		
			\addplot[color={black}, style={very thick}]
				table[x=x,y=gsn,col sep=comma,unbounded coords=jump]{./figs/data/sod_gamma.csv};
    		\addlegendentry{Exact $\gamma^*$ (Numerical)}
    
			\addplot[color={red!90!black}, style={ultra thick, dotted}]
				table[x=x,y=gsa,col sep=comma,unbounded coords=jump]{./figs/data/sod_gamma.csv};
    		\addlegendentry{Approximate $\gamma^{**}$ (Analytic)}
		\end{axis} 		
	\end{tikzpicture}}}
    \subfloat[Modified Sod]{
    \adjustbox{width=0.49\linewidth,valign=b}{    \begin{tikzpicture}[spy using outlines={rectangle, height=3cm,width=2.3cm, magnification=3, connect spies}]
		\begin{axis}[name=plot1,
		    axis line style={latex-latex},
		    axis x line=left,
            axis y line=left,
            clip mode=individual,
		    xlabel={$x$},
		    xtick={0,0.2,0.4,0.6,0.8,1},
    		xmin=0,
    		xmax=1,
    		x tick label style={
        		/pgf/number format/.cd,
            	fixed,
            	fixed zerofill,
            	precision=1,
        	    /tikz/.cd},
    		ylabel={$\gamma^*$, \textcolor{red!90!black}{$\gamma^{**}$}},
    		ytick={0.00, 0.01, 0.02, 0.03, 0.04, 0.05},
    		ymax=0.05,
    		y tick label style={
        		/pgf/number format/.cd,
        	    /tikz/.cd        	    },
    		legend style={at={(0.03,1)},anchor=north west,font=\small},
    		legend cell align={left},
    		style={font=\normalsize}]
    		
			\addplot[color={black}, style={very thick}]
				table[x=x,y=gsn,col sep=comma,unbounded coords=jump]{./figs/data/modsod_gamma.csv};
    
			\addplot[color={red!90!black}, style={ultra thick, dotted}]
				table[x=x,y=gsa,col sep=comma,unbounded coords=jump]{./figs/data/modsod_gamma.csv};
		\end{axis} 		
	\end{tikzpicture}}}
    \caption{\label{fig:sod_gamma} 
    Comparison of the mass correction factor $\gamma$ computed by the exact numerical approach ($\gamma^*$) and the approximate analytic form ($\gamma^{**}$) for the Sod shock tube (left) and modified Sod shock tube (right) at $t = 0.2$ using bounds preserving RK4 temporal integration, pseudospectral spatial discretization ($256$), and invariant domain preserving bounds.}
\end{figure}
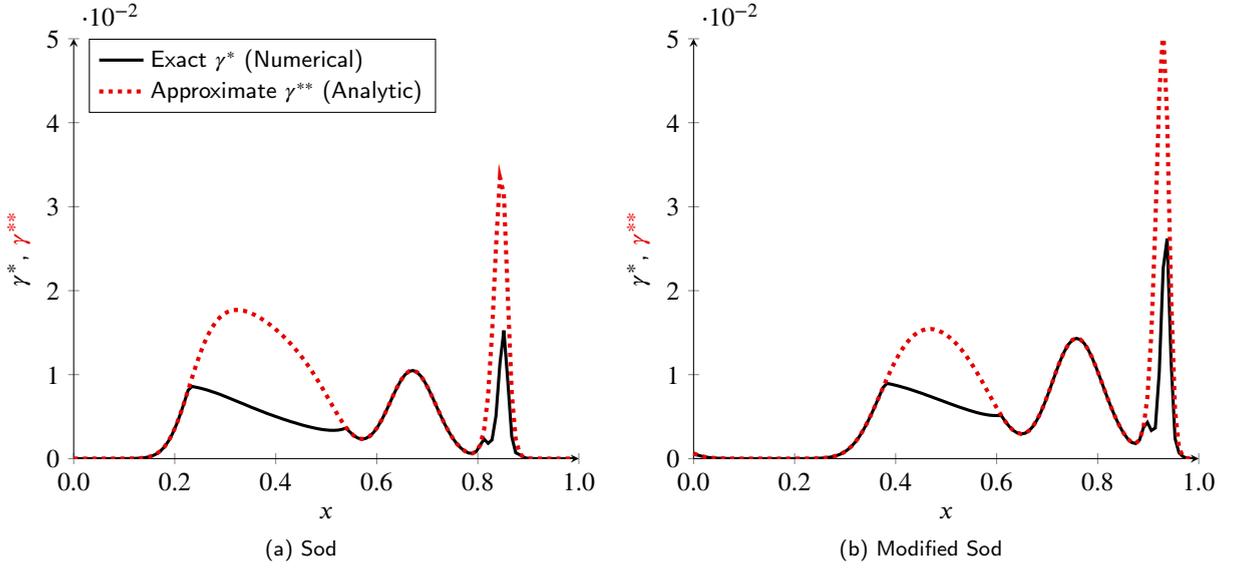

A comparison of the mass correction factor $\gamma$ computed by the exact numerical approach ($\gamma^*$) and the approximate analytic form ($\gamma^{**}$) is shown in \cref{fig:sod_gamma} at $t = 0.2$ for $N = 256$. The approximate form was generally in good agreement with the exact form in terms of behavior and magnitude, with excellent agreement in the region around the contact discontinuity. Additionally, the results confirm that $\gamma^*$ is bounded from above by $\gamma^{**}$. For either choice of computing the mass correction factor, the solution was essentially identical, and the use of the approximate analytic form did not violate the bounds at any point. 

\subsubsection{2D Riemann Problem}

For the extension to higher-dimensions, a two-dimensional Riemann problem was considered, introduced as case 12 in \citet{Liska2003}. The problem is defined on the domain $\Omega =[0,1]^2$ with the initial conditions given in \cref{fig:2dRP_ics}. Similarly to the Sod shock tube, the domain is periodized by reflecting the domain about the $x=0$ and $y=0$ axes. The contours of density computed using BP-RK4 temporal integration, pseudospectral spatial discretization ($N = 400^2$), and invariant domain preserving bounds with $\Delta t = 1{\cdot}10^{-4}$ are shown in \cref{fig:2dRP}. The results are in good agreement with the various methods in \citet{Liska2003}, with good resolution of the contact discontinuities and shock fronts and no observable spurious oscillations. A comparison with the results of a pseudospectral spatial discretization using an entropy viscosity approach from \citet{Guermond2011} is also shown in \cref{fig:2dRP}. Comparable results were seen with the BP-RK scheme even with a lower resolution and without an explicit shock-capturing approach. 

\begin{figure}[htbp!]
    \centering
    \adjustbox{width=0.45\linewidth,valign=b}{     \begin{tikzpicture}[spy using outlines={rectangle, height=3cm,width=2.3cm, magnification=3, connect spies}]
		\begin{axis}[name=plot1,
		    tick style={draw=none},
		    axis x line=left,
            axis y line=left,
            axis equal image,
            clip mode=individual,
    		xmin=0,
    		xmax=1,
    		xtick={0.0, 0.5, 1.0},
		    xlabel={$x$},
    		ymin=0,
    		ymax=1,
    		ytick={0.0, 0.5, 1.0},
		    ylabel={$y$},
    		ylabel style={rotate=-90},
    		style={font=\small},
    		scale = 1]

		    \draw[-] (axis cs:0.0, 0.0) -- (axis cs:1.0, 0.0);
		    \draw[-] (axis cs:1.0, 0.0) -- (axis cs:1.0, 1.0);
		    \draw[-] (axis cs:1.0, 1.0) -- (axis cs:0.0, 1.0);
		    \draw[-] (axis cs:0.0, 1.0) -- (axis cs:0.0, 0.0);
		    \draw[-] (axis cs:0.5, 0.0) -- (axis cs:0.5, 1.0);
		    \draw[-] (axis cs:0.0, 0.5) -- (axis cs:1.0, 0.5);
		    
		    ;
		  
		    \node at (axis cs:0.25,0.25) 
		    {\begin{tabular}{c} $\rho = 0.8$\\ $u = 0$ \\ $v = 0$ \\ $P = 1$ \end{tabular}};
		    \node at (axis cs:0.25,0.75) 
		    {\begin{tabular}{c} $\rho = 1.0$\\ $u = \frac{3}{\sqrt{17}}$ \\ $v = 0$ \\ $P = 1$ \end{tabular}};
		    \node at (axis cs:0.75,0.25) 
		    {\begin{tabular}{c} $\rho = 1.0$\\ $u = 0$ \\ $v = \frac{3}{\sqrt{17}}$ \\ $P = 1$ \end{tabular}};
		    \node at (axis cs:0.75,0.75) 
		    {\begin{tabular}{c} $\rho = \frac{17}{32}$\\ $u = 0$ \\ $v = 0$ \\ $P = 0.4$ \end{tabular}};
			    
		\end{axis}

	\end{tikzpicture}}
    \caption{\label{fig:2dRP_ics} Initial conditions for the 2D Riemann problem on the subdomain $\Omega = [0, 1]^2$. }
\end{figure}
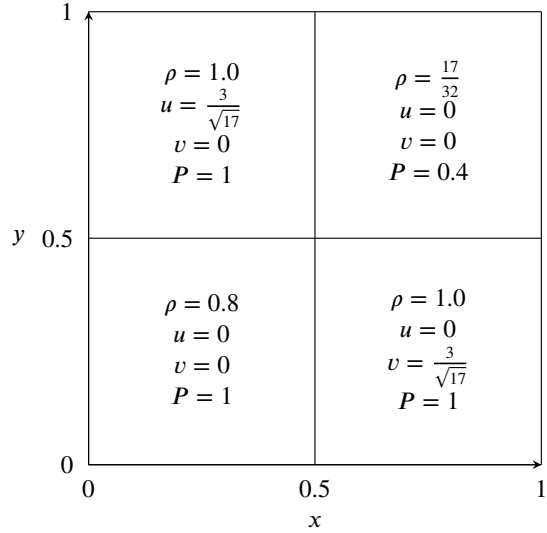
\begin{figure}[htbp!]
    \centering
    \subfloat[BP-RK4, $N = 400^2$]{\label{fig:2dRP_BPRK} \adjustbox{width=0.47\linewidth,valign=b}{\includegraphics[]{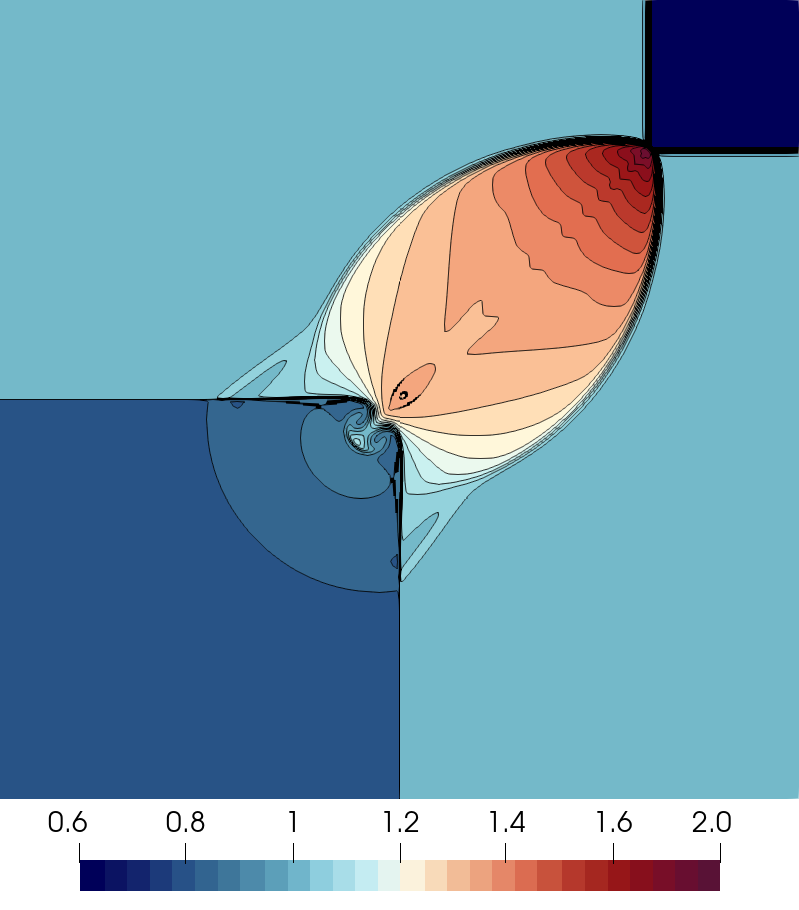}}}
    ~
    \subfloat[Entropy viscosity, $N = 600^2$ \citep{Guermond2011}]{\label{fig:2dRP_EV} \adjustbox{width=0.48\linewidth,valign=b}{\includegraphics[]{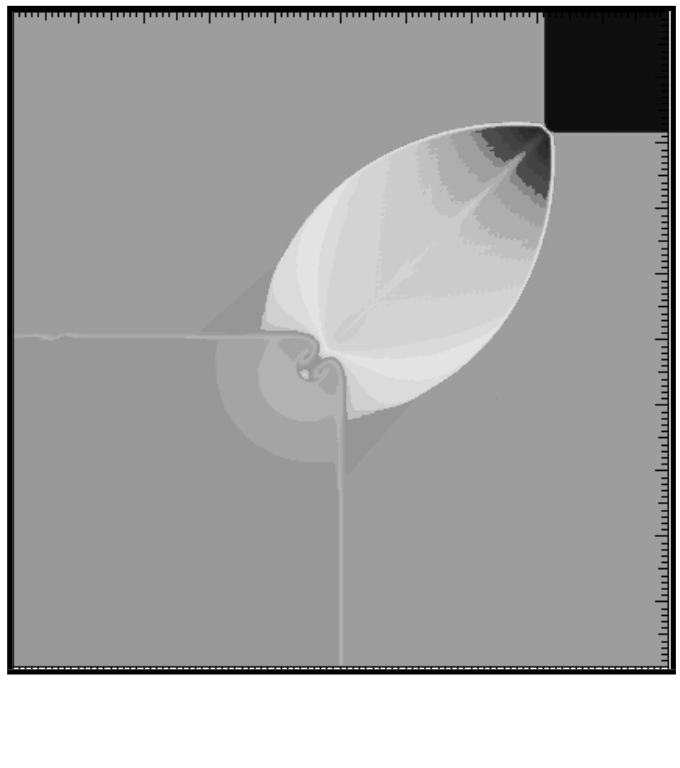}}}
    \caption{\label{fig:2dRP} Contours of density for the 2D Riemann problem at $t = 0.2$. \textit{Left}: Bounds preserving RK4 temporal integration, pseudospectral spatial discretization ($N = 400^2$), and invariant domain preserving bounds. \textit{Right}: Pseudospectral spatial discretization ($N = 600^2$) with entropy viscosity (\citet{Guermond2011}).}
\end{figure}

\subsubsection{Kelvin-Helmholtz Instability}
\begin{figure}[htb!]
    \centering
    \subfloat[$N = 100^2$]{\label{fig:kh_n100} \adjustbox{width=0.48\linewidth,valign=b}{\includegraphics[width=\textwidth]{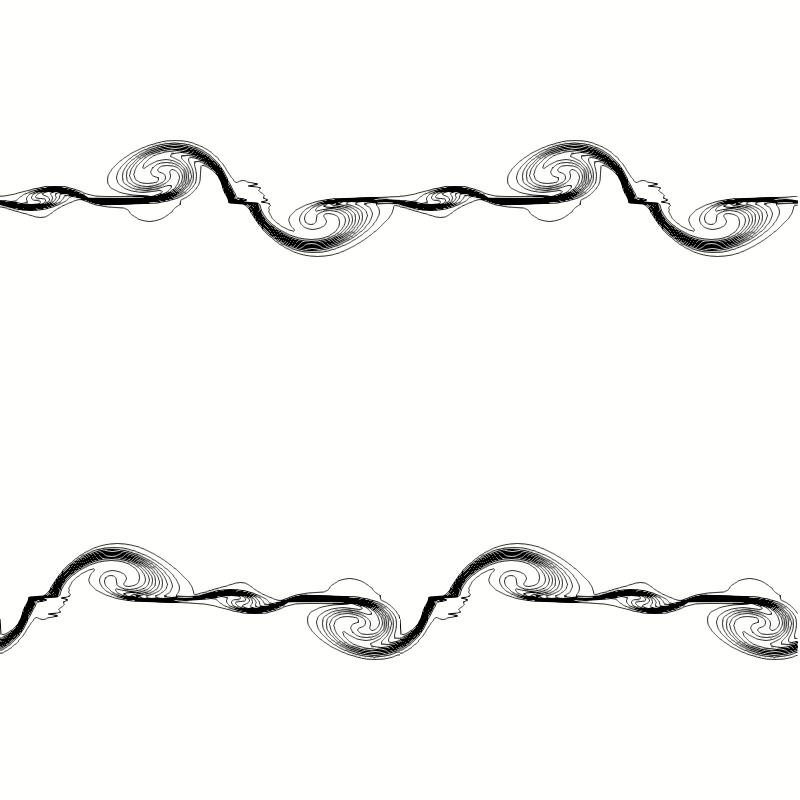}}}
    ~
    \subfloat[$N = 200^2$]{\label{fig:kh_n200} \adjustbox{width=0.48\linewidth,valign=b}{\includegraphics[width=\textwidth]{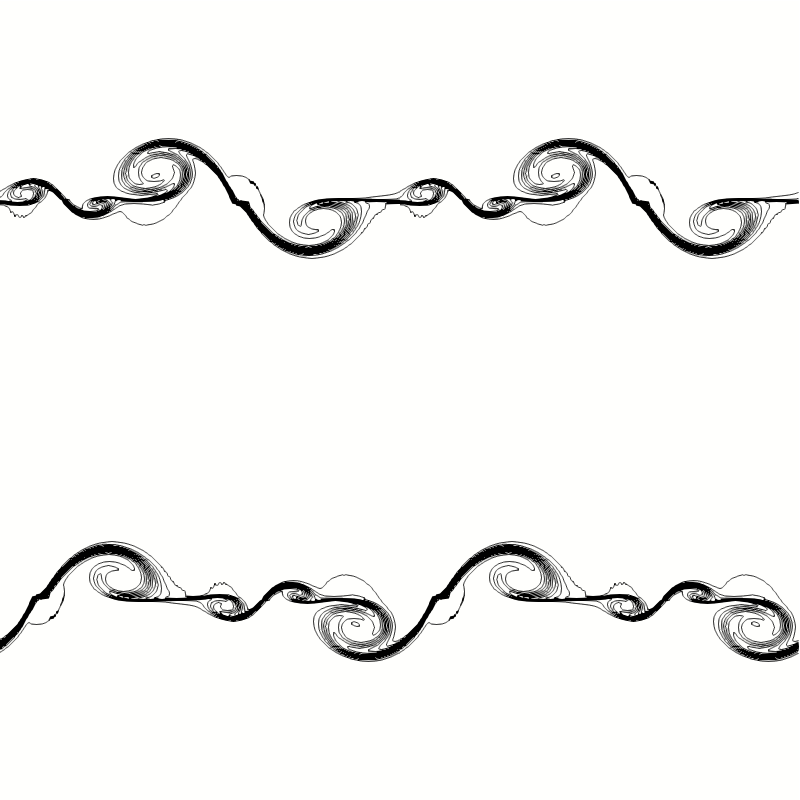}}}
    \newline
    \subfloat[$N = 300^2$]{\label{fig:kh_n300} \adjustbox{width=0.48\linewidth,valign=b}{\includegraphics[width=\textwidth]{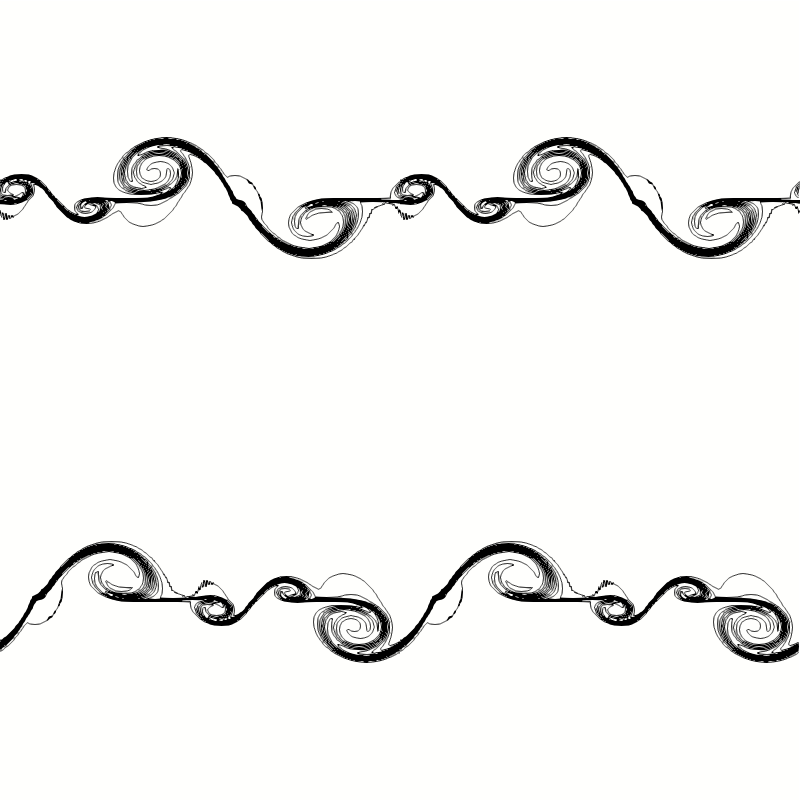}}}
    ~
    \subfloat[$N = 400^2$]{\label{fig:kh_n400} \adjustbox{width=0.48\linewidth,valign=b}{\includegraphics[width=\textwidth]{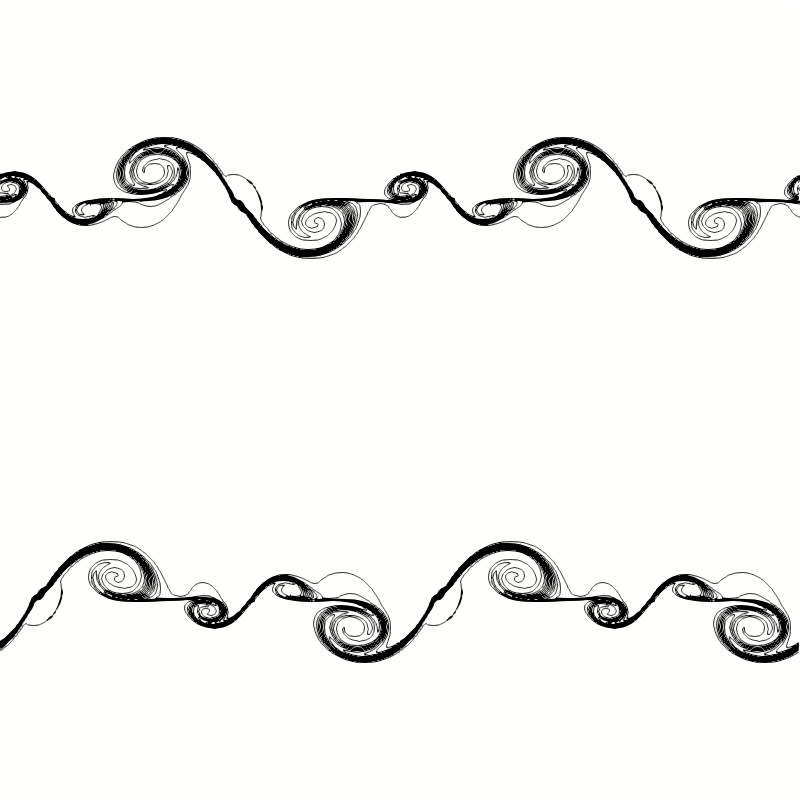}}}
    \newline
    \caption{\label{fig:kh} Isocontours of density for the Kelvin-Helmholtz instability at $t = 1$ using bounds preserving RK4 temporal integration, pseudospectral spatial discretization ($N = 100^2$, $200^2$, $300^2$, $400^2$), and invariant domain preserving bounds. }
\end{figure}

As a final assessment in the context of more complex flow physics, the vortical driven flow of a Kelvin-Helmholtz instability was considered. The problem is defined on the domain $\Omega =[-1,1]^2$ with the initial conditions given as

\begin{equation}
    \mathbf{q}_0(x) = \begin{cases}
        \mathbf{q}_l, &\mbox{if } |y + \phi(x)| \leq 0.5,\\
        \mathbf{q}_r, &\mbox{else},
    \end{cases} \quad \mathrm{given} \quad \mathbf{q}_l = \begin{bmatrix}
        2 \\ 0.5 \\ 0 \\ 2.5
    \end{bmatrix}, \quad \mathbf{q}_r = \begin{bmatrix}
        1 \\ -0.5 \\ 0 \\ 2.5
    \end{bmatrix},
\end{equation}
where $\phi(x)$ is an initial perturbation in the interface used to seed the instability. A sinusoidal perturbation of the form 
\begin{equation}
    \phi(x) = \alpha \sin \left (k\pi x \right)
\end{equation}
was used, where $\alpha = 10^{-2}$ is the amplitude and $k = 2$ is the frequency. The isocontours of density, shown as 10 equispaced isocontours on the range $[1,2]$, computed using BP-RK4 temporal integration, pseudospectral spatial discretization with $N = [100^2$, $200^2$, $300^2$, $400^2]$, and invariant domain preserving bounds with $\Delta t = [1{\cdot}10^{-3}$, $5{\cdot}10^{-4}$, $3{\cdot}10^{-4}$, $2{\cdot}10^{-4}]$ are shown in \cref{fig:kh}. The results show the rollup of vortices indicative of the Kelvin-Helmholtz instability, and the imposition of invariant domain preserving bounds did not excessively diffuse the smooth vortical structures in the flow while effectively suppressing spurious oscillations near density discontinuities. With progressively finer grids, the resolution of the vortices improved proportionally. 

\section{Conclusion}\label{sec:conclusion}
In this work, we introduced bounds preserving RK (BP-RK) schemes, a novel formulation of explicit RK temporal integration schemes that preserve any locally-defined quasiconvex set of bounds for the solution. These schemes operate on the basis of a nonlinear, bijective mapping between an admissible quasiconvex set of solutions and the real domain prior to temporal integration which is followed by an inverse mapping. The proposed techniques are generally applicable to a wide variety of problems, but the emphasis in this work was on nonlinear hyperbolic conservation laws, for which it was shown that an assortment of methods, such as positivity preserving, discrete maximum principle satisfying, entropy dissipative, and invariant  domain  preserving  schemes, could be recovered essentially independently of the spatial discretization. It was also shown, both analytically and experimentally, that the BP-RK schemes recover the order of accuracy of the underlying RK schemes upon which they are built and can be modified to preserve any linear invariant of the system. For many applications, the additional cost of the proposed schemes is almost negligible -- simply the evaluation of two mappings with closed-form solutions. To show the utility of these schemes, the results of the computation of nonlinear hyperbolic conservation laws with discontinuous solutions using a pseudospectral method without an explicit shock capturing approach were presented. Even though the spatial discretization scheme was ill-suited for discontinuous problems, the results of the BP-RK were on par with dedicated shock capturing schemes, showing good resolution of discontinuous features without any observable spurious oscillations. For general applications for which the spatial discretization schemes are better suited for the problems at hand, BP-RK schemes can potentially offer superior performance at a lower computational cost than dedicated spatial discretization schemes and allow for generalizability between various constraints without the need to derive and implement new discretizations. Future work may consider the application of the proposed scheme to implicit temporal integration which can present additional challenges in that the proper choice of bounds, particularly ones which rely on a local domain of influence, may become ambiguous and that the system may become ill-conditioned for solutions which approach the bounds as the magnitude of the mapping Jacobian would increase.

\section*{Acknowledgements}
\label{sec:ack}
This work was supported in part by the U.S. Air Force Office of Scientific Research via grant FA9550-21-1-0190 ("Enabling next-generation heterogeneous computing for massively parallel high-order compressible CFD") of the Defense University Research Instrumentation Program (DURIP) under the direction of Dr. Fariba Fahroo.

\bibliographystyle{unsrtnat}
\bibliography{references}



\end{document}